\documentclass[12pt]{article}
\usepackage{authblk}
\usepackage{amsthm,amsmath,amssymb}
\usepackage{cases}
\usepackage{dsfont}
\usepackage{multirow}
\usepackage{float}

\usepackage{dsfont}

\numberwithin{equation}{section}
\def\sumn{\sum^n_{i=1}}

\def\ep{\varepsilon}

\newtheorem{thm}{Theorem}[section]

\newtheorem{corollary}{Corollary}[section]

\newtheorem{remark}{Remark}[section]

\begin{document}

\date{}
\title{Bootstrapped Pivots for Sample and Population  Means and Distribution Functions  }

\author{Mikl\'{o}s Cs\"{o}rg\H{o}\thanks{mcsorgo@math.carleton.ca} and Masoud M. Nasari\thanks{mmnasari@math.carleton.ca}    \\
\small{School of Mathematics and Statistics of Carleton University}\\\small{ Ottawa, ON, Canada} }

\maketitle
\begin{abstract}
\noindent
In this paper we introduce the concept of \emph{bootstrapped pivots} for the sample and the  population means.
This is in  contrast to  the classical method of constructing bootstrapped confidence intervals for the population mean  via estimating the cutoff points  via  drawing a number of bootstrap sub-samples. We  show that this new method leads to constructing asymptotic   confidence intervals with  significantly smaller error in comparison to both of the traditional  $t$-intervals and  the classical bootstrapped confidence intervals. The approach taken in  this paper relates naturally to super-population modeling, as well as to estimating empirical and theoretical distributions.
\end{abstract}



\section{\normalsize{Introduction} }\label{section 1}
Unless stated otherwise, $X, X_{1},\ldots$ throughout   are assumed to be independent random variables with a common distribution function $F$ (i.i.d. random variables), mean $\mu:=E_X X$ and variance $0<\sigma^{2}:=E_X  (X-\mu)^2< +\infty$. Based on $X_1, \ldots,X_n$, a random sample on $X$,  for each integer  $n\geq 1$, define
\begin{equation*}
\bar{X}_{n}:=\sum_{i=1}^{n}X_i\big/n \ \textrm{and} \  S^{2}_{n}:=\sum_{i=1}^{n}(X_i -\bar{X}_{n})^2\big/ n,
\end{equation*}
the sample mean and sample variance, respectively, and consider the classical Student $t-$statistic
\begin{equation}\label{def. of T_n}
T_{n}(X):= \frac{\bar{X}_n}{S_n/\sqrt{n}}=\frac{\sum_{i=1}^n X_i}{S_n \sqrt{n}}
\end{equation}
that, in turn, on replacing $X_i$ by $X_i-\mu$, $1\leq i \leq n$, yields
\begin{equation}\label{equivalent to T_n}
T_{n}(X-\mu):= \frac{\bar{X}_n -\mu}{S_n/\sqrt{n}}=\frac{\sum_{i=1}^n (X_i-\mu)}{S_n \sqrt{n}},
\end{equation}
the classical Student pivot for the population mean $\mu$.
\par
Define now $T_{m_n}^{*}$ and  $G_{m_n}^{*}$, randomized versions of $T_{n}(X)$ and $T_{n}(X-\mu)$, respectively, as follows:

\begin{equation}\label{T^{*}}
T^{*}_{m_n}:= \frac{\sum_{i=1}^{n} \big( \frac{w^{(n)}_{i}}{m_{n}} -\frac{1}{n}   \big)  \ X_{i}   }{S_{n}\sqrt{\sum_{i=1}^{n} (\frac{w^{(n)}_i}{m_n}-\frac{1}{n})^{2}}  }=\frac{\sum_{i=1}^{n} \big( \frac{w^{(n)}_{i}}{m_{n}} -\frac{1}{n}   \big)  ( X_{i} -\mu)   }{S_{n}\sqrt{\sum_{i=1}^{n} (\frac{w^{(n)}_i}{m_n}-\frac{1}{n})^{2}}  },
\end{equation}

\begin{equation}\label{G^{*}}
G^{*}_{m_n}:= \frac{\sum_{i=1}^{n} \big| \frac{w^{(n)}_{i}}{m_{n}} -\frac{1}{n}   \big| \big( X_{i}-\mu \big)   }{S_{n}\sqrt{\sum_{i=1}^{n} (\frac{w^{(n)}_i}{m_n}-\frac{1}{n})^{2}}  },
\end{equation}
where,  the weights $(w^{(n)}_{1},\ldots,w^{(n)}_{n})$ have a multinomial distribution of size
$m_n:= \sum_{i=1}^n  w^{(n)}_{i}$ with respective probabilities  $1/n$, i.e.,
\begin{equation*}
  (w^{(n)}_{1},\ldots,w^{(n)}_{n})\  \substack{d\\=}\ \ multinomial(m_{n};\frac{1}{n},\ldots,\frac{1}{n}).
\end{equation*}
\par
The just introduced respective randomized $T_{m_n}^{*}$ and  $G_{m_n}^{*}$ versions of $T_{n}(X)$ and $T_{n}(X-\mu)$ can be computed via  re-sampling from  the set of indices $\{ 1,\ldots,n\}$ of   $X_1,\ldots,X_n$ with replacement $m_n$ times so  that, for each $1\leq i \leq n$,  $w^{(n)}_i$ is  the count of the  number of times the index $i$ of $X_i$ is chosen in this  re-sampling process. 

\begin{remark}\label{Remark 1}
In view of the preceding definition of $w^{(n)}_{i}$, $1\leq i \leq n$, they form a row-wise independent triangular array of random variables such that  $\sum_{i=1}^n  w^{(n)}_{i}=m_n$ and, for each $n\geq 1$,
$$(w^{(n)}_{1},\ldots,w^{(n)}_{n})\  \substack{d\\=}\ \ multinomial(m_{n};\frac{1}{n},\ldots,\frac{1}{n}),
$$
i.e., the weights  have a multinomial distribution of size $m_n$ with respective probabilities $1/n$.
Clearly,  for each $n$, $w^{(n)}_{i}$    are independent from the random sample $X_{i}$, $1\leq i \leq n$.
Weights denoted by $w_{i}^{(n)}$ will stand for triangular multinomial random variables in this context throughout.
\end{remark}

\par
We note that our approach to the bootstrap is to benefit  (cf. Corollary \ref{The Rate} and Remark \ref{Remark 2.1}) from the refinement provided by the respective  scalings  $\big| \frac{w^{(n)}_{i}}{m_{n}} -\frac{1}{n}   \big|\big/ (S_{n}\sqrt{\sum_{i=1}^{n} (\frac{w^{(n)}_i}{m_n}-\frac{1}{n})^{2}}  ) $ of $(X_i -\mu)$'s and $ (\frac{w^{(n)}_{i}}{m_{n}} -\frac{1}{n} )  \big/ (S_{n}\sqrt{\sum_{i=1}^{n} (\frac{w^{(n)}_i}{m_n}-\frac{1}{n})^{2}}  ) $  of $X_i$'s  in $G^{*}_{m_n}$ and $T^{*}_{m_n}$ which result from the re-sampling, as compared to  $1/(S_n \sqrt{n})$ of the classical central limit theorem (CLT) simultaneously for $T_n (X)$  and $T_n (X-\mu)$.  Viewed as  pivots, $G_{m_n}^{*}$ and $T_{m_n}^{*}$ are directly related to the respective parameters  of interest, i.e.,    the population mean $\mu=E_X X$ and   the  sample mean $\bar{X}_n$,   when, e.g.,  it is to be estimated  as
the sample mean of an imaginary random sample from an infinite super-population. This approach  differs fundamentally from the classical use of the bootstrap in which re-sampling is used to capture the sampling distribution of the pivot $T_{n} (X-\mu)$ for $\mu$ using $T^{*}_{m_n}$ that, by definition, is no longer  related to $\mu$ (cf. Section \ref{section pivot for bar{X}}).

\par
In the literature of the bootstrap,  the asymptotic behavior of $T_{m_n}^{*}$ is usually studied  by  conditioning on the observations on $X$. In this paper  we  use the method of first conditioning on the weights $w_{i}^{(n)}$ for both $T_{m_n}^{*}$ and $G_{m_n}^{*}$, and then conclude also their limit law in terms  of the joint distribution of the observations and the weights. We use the same approach for studying the asymptotic laws of the sample and population distribution functions.

 For similarities and differences between  the two methods of      conditioning, i.e.,  on the observations and on the weights,   for studying   the convergence in  distribution of  $T_{m_n}^{*}$,  we refer to Cs\"{o}rg\H{o} \emph{et al}. \cite{Scand}.

\textbf{Notations}. Let $(\Omega_X,\mathfrak{F}_X,P_X)$ denote  the probability space of the random variables  $X,X_1,\ldots$, and $(\Omega_w,\mathfrak{F}_w,P_w)$ be  the probability space on which the weights $\big(w^{(1)}_1,(w^{(2)}_1,w^{(2)}_{2}),\ldots,(w^{(n)}_1,\ldots,w^{(n)}_{n}),\ldots \big)$ are defined. In view of the independence of these two sets of random variables, jointly they live on the direct product probability space $(\Omega_X \times \Omega_w, \mathfrak{F}_X \otimes \mathfrak{F}_w, P_{X,w}=P_X .\ P_w )$. For each $n\geq 1$, we also let  $P_{.|w}(.)$ and $P_{.|X}(.)$  stand for the conditional probabilities  given $\mathfrak{F}^{(n)}_w:=\sigma(w^{(n)}_1,\ldots,w^{(n)}_{n})$ and  $\mathfrak{F}^{(n)}_X:=\sigma(X_1,\ldots,X_n)$, respectively, with corresponding conditional expected values $E_{.|w}(.)$ and $E_{.|X}(.)$.
\par
We are to  outline  now our view of $T^{*}_{m_n}$ and $G^{*}_{m_n}$ in terms of their related roles in statistical inference as they are studied in this exposition.
\par
To begin with, $T^{*}_{m_n}$ is one of the well studied bootstrapped versions of the Student $t$-statistic $T_{n}(X)$ for constructing bootstrap   confidence intervals for $\mu$. Unlike $T_{n}(X)$ however, that can be transformed into $T_{n}(X-\mu)$, the Student pivot for $\mu$ as in (\ref{equivalent to T_n}), its bootstrapped version $T^{*}_{m_n}$ does not have this straightforward  property, i.e., it does not yield a pivotal quantity for the population mean $\mu=E_{X} X$ by simply replacing each  $X_i$ by $X_i-\mu$ in its definition. Indeed, it takes considerably more effort to construct and establish the validity of bootstrapped $t-$intervals for $\mu$ by means of approximating the sampling distribution of $T_{n}(X-\mu)$ via that of $T_{m_n}^{*}$. For references and recent contributions  along these lines we refer to Cs\"{o}rg\H{o} \emph{et al.} \cite{Scand} and to Section \ref{Comparison section} of this exposition.

 \par
Along with other results, Section \ref{berry-esseen section} contains the main objective of this paper  which is to show  that when $E_{X}|X|^r <+\infty$, $r\geq 3$, then the confidence bound for the population mean $\mu$,   $G^{*}_{m_n}\leq z_{1-\alpha}$, where $z_{1-\alpha}$ is the $(1-\alpha)$th percentile of a standard normal distribution,  converges  to its nominal probability  coverage $1-\alpha$ at a significantly  better rate of at most  $O(n^{-1})$ (cf. Theorem \ref{Berry-Esseen},  Corollary  \ref{The Rate} and Remark \ref{Remark 2.1})  than that of the traditional  $t-$interval $T_{n}(X-\mu) \leq z_{1-\alpha}$ which is of order $O(n^{-1/2})$. This improvement results from   incorporating   the additional randomness provided by the weights  $w_{i}^{(n)}$. In fact, it is shown in Section \ref{berry-esseen section}
   that, when $E_{X}|X|^r <+\infty$, $r\geq 3$,  the conditional distribution functions, given $w_{i}^{(n)}$'s,   of   both $G^{*}_{m_n}$ and $T^{*}_{m_n}$ converge   to the standard normal distribution  at a rate that is  at best $O(n^{-1})$. A numerical study of the  results of  Section \ref{berry-esseen section} is also presented there.

\par
In Section \ref{section pivot for mu}, $G^{*}_{m_n}$ is introduced as a natural direct  asymptotic pivot for the population mean $\mu=E_{X} X$,  while in Section \ref{section pivot for bar{X}}, $T^{*}_{m_n}$  is studied on its own as a natural pivot for the sample mean $\bar{X}_n$.      Furthermore, in Section \ref{section super population}, both  $G^{*}_{m_n}$ and $T^{*}_{m_n}$ will be seen to relate  naturally to sampling from a finite population that is viewed as an imaginary random sample of size $N$ from an infinite super-population (cf., e.g.,  Hartley and Sielken Jr. \cite{Hartley}), rather than to sampling from a finite population of real numbers of size $N$.

\par

Confidence intervals (C.I.'s)  are established in a similar vein for the empirical and theoretical distributions in Section \ref{CLT empirical}.

\par
In Section \ref{Comparison section},  we compare the   rate at which the actual probability  coverage of   the C.I. $G^{*}_{m_n}\leq z_{1-\alpha}$ converges  to its nominal level, to the  classical bootstrap  C.I.  for $\mu$,  in which $T_{n}(X-\mu)$ is the pivot and the cutoff points  are estimated by the means of $T^{*}_{m_n}(1),\ldots,T^{*}_{m_n}(B)$ which are computed  versions of $T^{*}_{m_n}$, based on  $B$ independent bootstrap sub-samples of size $m_n$. In this paper we specify the phrase ``classical bootstrap C.I." for the aforementioned  C.I.'s based on bootstrapped cutoff points.

   In Section \ref{Comparison section} we also introduce a more accurate version of the classical bootstrap C.I. that is  based on
   a fixed number of bootstrap sub-samples $B$. In other words,  for our version of the classical   bootstrap C.I. to  converge  to its nominal probability  coverage,  $B$ does not have to be particularly large.  This result coincides with one of the results obtained by   Hall \cite{Hall B} in this regard.
   Also, as it will be seen in Section \ref{numerical comarison of the three}, when $B=9$ bootstrap sub-samples are drawn, in our   version of the classical bootstrap,  the \emph{precise} nominal coverage probability 0.9000169      coincides with the one of level 0.9  proposed  in Hall \cite{Hall B} for the same number of bootstrap replications $B=9$.
   Our investigation of the classical  bootstrap C.I.'s for $\mu$ in Section \ref{Comparison section} shows that the rate at which their actual probability coverage  approaches their  nominal one  is no better than that of the classical CLT, i.e., $O(n^{-1/2})$.

\par
Section \ref{numerical comarison of the three} is devoted to presenting a numerical comparison between the    three methods of constructing  C.I.'s for the population mean $\mu=E_X X$,  namely  the C.I. based on our bootstrapped  pivot $G^{*}_{m_n}$, the traditional  C.I. based on the pivot $T_{n}(X-\mu)$ and the classical bootstrap C.I..     As it will be seen in Section \ref{numerical comarison of the three},   the use of the bootstrapped pivot $G^{*}_{m_n}$ tends to significantly  outperform the other two methods by   generating  values closer to the nominal coverage probability more often.

\par

 The proofs are given in Section \ref{Proofs}.

\section{\normalsize{The rate of convergence of the CLT's for   $G_{m_n}^{*}$ and $T_{m_n}^{*}$ }}\label{berry-esseen section}
The asymptotic behavior  of the bootstrapped $t$-statistics has been extensively studied in the literature.  A  common feature of  majority of the  studies that  precede paper \cite{Scand} is that they were conducted by means of conditioning on a   given  random sample on $X$. In Cs\"{o}rg\H{o} \emph{et al.} \cite{Scand}, the asymptotic normality of  $T^{*}_{m_n}$,  when $0< \sigma^{2}:=E_{X}(X -\mu)^2< +\infty$, was investigated  via introducing the  approach  of conditioning on the bootstrap weights $w^{(n)}_{i}$'s (cf. (\ref{(1.11)}) in this exposition),  as well as via revisiting the method of  conditioning on the sample when $X$ is in the domain of attraction of the normal law, possibly with $E_{X} X^2 =+\infty$. Conditioning on the weights,  the asymptotic normality of the direct pivot $G^{*}_{m_n}$ for the population mean $\mu$,  when $0<\sigma^2 < +\infty$,  can be proven by simply replacing $(\frac{w^{(n)}_i}{m_n}-\frac{1}{n})$ by $\big|\frac{w^{(n)}_i}{m_n}-\frac{1}{n}\big|$ in the proof of part (a) of Corollary 2.1 in  paper \cite{Scand}, leading to concluding (\ref{(1.23)}) as $n,m_n \rightarrow +\infty$ so that $m_n=o(n^2)$.

 \par
One efficient tool to control  the error when   approximating the distribution function  of a statistic with that  of a standard normal random variable  is provided by   Berry-Ess\'{e}en type inequalities  (cf. for example Serfling \cite{Berry-Esseen}),
 which provides  an upper bound  for the error  for any finite number of observations
in hand. It is well known that, as  the sample size  $n$ increases to infinity,   the rate at which the Berry-Esse\'{e}n    upper bound vanishes is $n^{-1/2}$.

\par
Our Berry-Esse\'{e}n type inequality  for the respective conditional, given $w^{(n)}_i$'s, distributions of  $G^{*}_{m_n}$ and $T^{*}_{m_n}$,  as in (\ref{T^{*}}) and (\ref{G^{*}}) respectively, reads as follows.

\begin{thm}\label{Berry-Esseen}
Assume that $E_{X}|X|^{3}<+\infty$ and let $\Phi(.)$ be the standard normal  distribution function.
Also, for arbitrary positive numbers $\delta,\varepsilon$, let $\varepsilon_1, \varepsilon_2>0$ be so that $\delta> (\varepsilon_1/\varepsilon)^2+P_{X}(|S^{2}_n -\sigma^2|>\varepsilon^{2}_1 )+\varepsilon_2 >0$, where, for $t\in \mathbb{R}$, $\Phi(t-\varepsilon)-\Phi(t)>- \varepsilon_2$ and $\Phi(t+\varepsilon)-\Phi(t)<  \varepsilon_2$.
 Then, for all $n,m_n$    we have

\begin{eqnarray*}
(A)&&P_{w}\big\{ \sup_{-\infty < t < +\infty}\Big| P_{X|w} (G^{*}_{m_n}\leq t)-\Phi(t)  \Big|>\delta \big\}~~~~~~~~~~~~~~~~~~~~~~~~~~~~\\
&&\leq \delta^{-2}_n(1-\varepsilon)^{-3}(1-\frac{1}{n})^{-3 } (\frac{n}{m^{3}_n}+\frac{n^{2}}{m^{3}_n})\{ \frac{15 m^{3}_n}{n^3}+ \frac{25 m^{2}_n}{n^2}+\frac{m_n}{n} \}  \\
&&+ \varepsilon^{-2}\frac{m^{2}_n}{ (1-\frac{1}{n})} \Big\{ \frac{1-\frac{1}{n}}{n^3 m^{3}_n } + \frac{(1-\frac{1}{n})^4}{m^{3}_n}  + \frac{(m_n -1)(1-\frac{1}{n})^2}{n m^{3}_n} +  \frac{4(n-1)}{n^3 m_n}       +\frac{1}{m^{2}_n} \\
&&- \frac{1}{n m^{2}_n} + \frac{n-1}{n^{3} m^{3}_{n}} + \frac{4(n-1)}{n^2 m^{3}_{n}} - \frac{(1-\frac{1}{n})^2}{m^{2}_n}        \Big\},
\end{eqnarray*}
and also
\begin{eqnarray*}
(B)&&P_{w}\big\{ \sup_{-\infty < t < +\infty}\Big| P_{X|w} (T^{*}_{m_n}\leq t)-\Phi(t)  \Big|>\varepsilon \big\}~~~~~~~~~~~~~~~~~~~~~~~~~~~~\\
&&\leq \delta^{-2}_n(1-\varepsilon)^{-3}(1-\frac{1}{n})^{-3 } (\frac{n}{m^{3}_n}+\frac{n^{2}}{m^{3}_n})\{ \frac{15 m^{3}_n}{n^3}+ \frac{25 m^{2}_n}{n^2}+\frac{m_n}{n} \}  \\
&&+ \varepsilon^{-2}\frac{m^{2}_n}{ (1-\frac{1}{n})} \Big\{ \frac{1-\frac{1}{n}}{n^3 m^{3}_n } + \frac{(1-\frac{1}{n})^4}{m^{3}_n}  + \frac{(m_n -1)(1-\frac{1}{n})^2}{n m^{3}_n} +  \frac{4(n-1)}{n^3 m_n}       +\frac{1}{m^{2}_n} \\
&&- \frac{1}{n m^{2}_n} + \frac{n-1}{n^{3} m^{3}_{n}} + \frac{4(n-1)}{n^2 m^{3}_{n}} - \frac{(1-\frac{1}{n})^2}{m^{2}_n}        \Big\},
\end{eqnarray*}
where
\begin{equation*}
\delta_n:= \frac{\delta- (\varepsilon_1/\varepsilon)^2 -P_{X} ( |S^{2}_n -\sigma^2|>\varepsilon^{2}_1 )+\varepsilon_2   }{C E_{X}|X-\mu|^{3}/\sigma^{3/2}},
\end{equation*}
with $C$ being  the Berry-Esse\'{e}n universal constant.
\end{thm}

\par

The following result, a corollary to Theorem \ref{Berry-Esseen}, gives the rate of convergence of the respective conditional  CLT's for $G^{*}_{m_n}$ and $T^{*}_{m_n}$.

\begin{corollary}\label{The Rate}
Assume that $E_{X}|X|^{r}<+\infty$, $r\geq 3$. If $n, \ m_n \rightarrow +\infty$ in such a way that $m_n=o(n^2)$,  then, for  arbitrary $\delta>0$,  we have
\begin{equation*}
(A) \ P_{w}\big\{ \sup_{-\infty < t < +\infty}\Big| P_{X|w} (G^{*}_{m_n}\leq t)-\Phi(t)  \Big|>\delta \big\} = O\Big( \max\{ \frac{m_n}{n^2}, \frac{1}{m_n} \}  \Big), ~~~~~~~~~~~~~
\end{equation*}
\begin{equation*}
(B) \ P_{w}\big\{ \sup_{-\infty < t < +\infty}\Big| P_{X|w} (T^{*}_{m_n}\leq t)-\Phi(t)  \Big|>\delta \big\} = O\Big( \max\{ \frac{m_n}{n^2}, \frac{1}{m_n}  \}  \Big). ~~~~~~~~~~~~~
\end{equation*}
Moreover, if $E_{X} X^4 <+\infty$, if $n,m_n \to +\infty$ in such a way that $m_n=o(n^2)$ and $n=o(m^{2}_n)$ then, for  $\delta>0$,  we also have
\begin{equation*}
(C) \ P_{w}\big\{ \sup_{-\infty < t < +\infty}\Big| P_{X|w} (G^{**}_{m_n}\leq t)-\Phi(t)  \Big|>\delta \big\} = O\Big( \max\{ \frac{m_n}{n^2}, \frac{1}{m_n},\frac{n}{m^{2}_{n}} \}  \Big), ~~~~~~~~~~~~~
\end{equation*}
\begin{equation*}
(D) \ P_{w}\big\{ \sup_{-\infty < t < +\infty}\Big| P_{X|w} (T^{**}_{m_n}\leq t)-\Phi(t)  \Big|>\delta \big\} = O\Big( \max\{ \frac{m_n}{n^2}, \frac{1}{m_n},\frac{n}{m^{2}_{n}} \}  \Big), ~~~~~~~~~~~~~
\end{equation*}
where
\begin{eqnarray}
G^{**}_{m_n}&:=& \frac{\sum_{i=1}^{n} \big| \frac{w^{(n)}_{i}}{m_{n}} -\frac{1}{n}   \big| \big( X_{i}-\mu \big)   }{S^{*}_{m_n}\sqrt{\sum_{i=1}^{n} (\frac{w^{(n)}_i}{m_n}-\frac{1}{n})^{2}}  }, \label{G^{**}} \\
T^{**}_{m_n}&:=& \frac{\bar{X}^{*}_{m_n}-\bar{X}_n   }{S^{*}_{m_n}\sqrt{\sum_{i=1}^{n} (\frac{w^{(n)}_i}{m_n}-\frac{1}{n})^{2}}  }=\frac{\sum_{j=1}^{n} \big( \frac{w^{(n)}_{i}}{m_{n}} -\frac{1}{n}   \big)   X_{i}   }{S^{*}_{m_n}\sqrt{\sum_{i=1}^{n} (\frac{w^{(n)}_i}{m_n}-\frac{1}{n})^{2}}  } \label{T^{**}}
\end{eqnarray}
and $S_{m_n}^{*^2}$ is the bootstrapped sample variance  defined as
\begin{equation}\label{def. of S^*}
S^{*^{2}}_{m_n}:= \sum_{i=1}^n w_{i}^{(n)} \big( X_i- \bar{X}^{*}_{m_n} \big)^2 \big/ m_n,
\end{equation}
where
\begin{equation}\label{bootstrap mean}
\bar{X}^{*}_{m_n}:=\sum_{i=1}^n w^{(n)}_i X_i/m_{n}.
\end{equation}
\end{corollary}

\begin{remark}\label{Remark 2.1}
On taking $m_n=n$, when $E_{X}|X|^{3}<+\infty$, the rates of convergence of Corollary \ref{The Rate} for both $G^{*}_{m_n}$ and $T^{*}_{m_n}$ are of order
$O( n^{-1})$.  The same is true for $G^{**}_{m_n}$ and $T^{**}_{m_n}$ for $m_n=n$ when $E_{X} X^{4}<+\infty$.
\end{remark}

\begin{remark}\label{Remark 2.2}

When $E_{X} X^4 <+\infty$, the extra term of  $n/m_{n}^{2}$ which appears in the rate of convergence of   $G^{**}_{m_n}$ and $T^{**}_{m_n}$ in (C) and (D) of Corollary \ref{The Rate}, is the rate at which $P_{w}\big\{ P_{X|w}\big( |S^{*^2}_{m_n}-S^{2}_{n}|>\varepsilon_1\big) > \varepsilon_2 \big\}$ approaches zero as $n,m_n \to +\infty$,  where  $\varepsilon_1$ and $\varepsilon_2$ are arbitrary positive numbers.
\end{remark}

\par
Furthermore, in view of Lemma 1.2 of S. Cs\"{o}rg\H{o} and Rosalsky \cite{Csorgo and Rosalsky},  the conditional CLT's resulting from (A), (B), (C) and (D) of Corollary \ref{The Rate}  imply respective unconditional CLT's in terms of $P_{X,w}$. Moreover, when $m_n=n\to +\infty$, we conclude

\begin{eqnarray}
\sup_{-\infty < t < +\infty}  \big|  P_{X,w} (G^{*}_{m_n}\leq t)-\Phi(t) \big| &=& O(1/n), \label{joint dis1}\\
\sup_{-\infty < t < +\infty}\Big| P_{X,w} (T^{*}_{m_n}\leq t)-\Phi(t)  \Big| &=& O(1/n), \label{joint dis2}\\
\sup_{-\infty < t < +\infty}  \big|  P_{X,w} (G^{**}_{m_n}\leq t)-\Phi(t) \big| &=& O(1/n), \label{joint dis3}\\
\sup_{-\infty < t < +\infty}  \big|  P_{X,w} (T^{**}_{m_n}\leq t)-\Phi(t) \big| &=& O(1/n), \label{joint dis4}
\end{eqnarray}
where, (\ref{joint dis1}) and (\ref{joint dis2}) hold true when $E_{X} |X|^3 <+\infty$ and (\ref{joint dis3}) and (\ref{joint dis4}) hold true when $E_{X} X^4 <+\infty$.

\par

In the following Table 1 we use the software R to  present some numerical illustrations of our Theorem \ref{Berry-Esseen} and its Corollary \ref{The Rate} for $G^{*}_{m_n}$ as compared to the CLT  for   $T_n(X-\mu)$.

The conditional  probability of $G^{*}_{m_n}$, with $m_n=n$,  given $w_{i}^{(n)}$'s, is approximated by its empirical counterpart. On  taking $m_n=n$, the procedure for each $n$   is  generating  a realization of the weights $(w_1^{(n)},\ldots,w_{n}^{(n)})$ first,  and then  generating  500 sets of $ X_1,\ldots,X_n$ to compute the conditional  empirical distribution $\sum_{t=1}^{500} \mathds{1}(G^{*}_{m_n}\leq 1.644854)/500$. Simultaneously, and for the same generated 500 sets of $ X_1,\ldots,X_n$, we also compute the empirical distribution $\sum_{t=1}^{500} \mathds{1}(T_{n}(X-\mu) \leq 1.644854)/500$, where, here and throughout this paper, $\mathds{1}(.)$ stands for the indicator function.
 \par
  This cycle is repeated 500 times, i.e., for 500 realizations  $(w^{(n)^{(s)}}_1,\ldots,w^{(n)^{(s)}}_n)$, $1\leq s \leq 500$. For each $s$, 500 sets   $X_{1}^{(s,t)},\ldots,X_{n}^{(s,t)}$, $1\leq t \leq 500$ are generated. Then we compute the following empirical  values.

\begin{eqnarray*}\nonumber
emp_{X|w}G^{*}&:=&  \frac{ \sum_{s=1}^{500}
\mathds{1} \big( \big| \frac{\sum_{t=1}^{500} \mathds{1}(G^{*^{(s,t)}}_{m_n}\leq 1.644854)}{500}- 0.95   \big| \leq 0.01 \big) } {500},\\
emp_{X}T_n(X-\mu)&:=& \frac{ \sum_{s=1}^{500} \mathds{1}\big( \big| \frac{\sum_{t=1}^{500} \mathds{1}(T^{(s,t)}_{n}(X-\mu)\leq 1.644854)} {500}- 0.95  \big| \leq 0.01 \big) } {500},
\end{eqnarray*}
where, $G^{*^{(s,t)}}_{m_n}$ stands for the computed value of $G^{*}_{m_n}$, with $m_n=n$, based on
 $(w^{(n)^{(s)}}_1,\ldots,w^{(n)^{(s)}}_n)$ and $X_{1}^{(s,t)},\ldots,X_{n}^{(s,t)}$ and  $T^{(s,t)}_{n}(X-\mu)$ stands for the computed value of $T_{n}(X-\mu)$ based on $X_{1}^{(s,t)},\ldots,X_{n}^{(s,t)}$.
\par
The numerical  results, which  are presented in the following table,  indicate a  significantly better performance of $G^{*}_{m_n}$, with $m_n=n$, in comparison to that of the Student  $t$-statistic $T_n(X-\mu)$.
\\
\begin{table}[ht]\label{Table 1}
\caption{Comparing $G^{*}_{m_n}$, with $m_n=n$, to $T_{n}(X-\mu)$ }
\begin{center}
\begin{tabular}{ c|c|c|c  }
\hline \hline
  Distribution &  $n$ & $emp_{X|w}G^{*}$ & $emp_{X}T_n(X-\mu)$  \\
\hline
  \multirow{3}{*}{$Poisson(1)$} &20& 0.552 &0.322 \\
                                &30& 0.554 &0.376   \\
                                &40& 0.560 &0.364 \\
                         \hline

 \multirow{3}{*}{$Lognormal(0,1)$}     &20& 0.142& 0.000  \\
                                       &30& 0.168& 0.000   \\
                                       &40& 0.196& 0.000 \\
                         \hline

 \multirow{3}{*}{$Exponentia(1)$}   &20&  0.308 & 0.016  \\
                                    &30&  0.338& 0.020  \\
                                    &50&  0.470 & 0.094   \\
\hline
\end{tabular}
\end{center}
\end{table}

\par
In view of (\ref{joint dis1})-(\ref{joint dis4}), we also compare the rate of convergence  of actual coverage probability of  the confidence bounds   $G_{m_n}^{*}\leq z_{1-\alpha}$, with $m_n=n$ and $0< \alpha<1$, to those of the traditional confidence bounds   $T_n (X-\mu)\leq z_{1-\alpha}$ and  the classical bootstrap  confidence bounds of size $1-\alpha$ (cf. Table 2 in Section \ref{numerical comarison of the three}).

\section{\normalsize{Bootstrapped asymptotic pivots for the population mean $\mu$}}\label{section pivot for mu}
We are now to present $G_{m_n}^{*}$ of (\ref{G^{*}}), and some further versions of it, as \emph{direct} asymptotic bootstrap pivots for the population mean $\mu=E_{X} X$ when $0<\sigma^2:=E_{X}(X-\mu)^2<+\infty$.
\par
We note that for the numerator term of $G_{m_n}^{*}$ we have
$$E_{X|w}\big( \sum_{i=1}^n |\frac{w_{i}^{(n)}}{m_n}-\frac{1}{n}|(X_i-\mu) \big)=0.$$
Furthermore, given $w^{(n)}_i$'s, for the bootstrapped weighted average
\begin{equation}\label{(1.21)}
{\sum_{i=1}^n |\frac{w_{i}^{(n)}}{m_n}-\frac{1}{n}| (X_i-\mu)}=: \bar{X}_{m_n}^{*}(\mu),
\end{equation}
mutatis mutandis in verifying (\ref{(5.30)}) in Appendix 1, we conclude that  when the original sample size $n$ \emph{is fixed} and $m:=m_n$, then,  \emph{as} $m \to +\infty$, we have 
\begin{equation}\label{(1.22)}
\bar{X}_{m_n}^{*}(\mu)=\bar{X}_{m}^{*}(\mu) \to 0 \ in \ probability-P_{X,w},
\end{equation}
and the same holds true if $n\to +\infty$ as well.
\par
In view of (\ref{(1.21)})
\begin{equation}\label{(1.22)'}
\frac{\sum_{i=1}^n |\frac{w_{i}^{(n)}}{m_n}-\frac{1}{n}| X_i}{\sum_{j=1}^n |\frac{w_{j}^{(n)}}{m_n}-\frac{1}{n}|}=: \bar{X}_{n,m_n}^{*}.
\end{equation}
is an unbiased estimator for $\mu$ with respect to $P_{X|w}$.

\par
It can be shown that when $E_X X^2<+\infty$, as  $n,m_n \to +\infty$ such that $m_n=o(n^2)$, $\bar{X}_{n,m_n}^{*}$ is a consistent estimator for the population mean $\mu$ in terms of $P_{X,w}$, i.e.,
\begin{equation}\label{doubleprime}
\bar{X}_{n,m_n}^{*} \to \mu \ in \ probability-P_{X,w}.
\end{equation}
In Appendix 1 we give a  direct proof for  (\ref{doubleprime})  for the important case when $m_n=n$, for which the  CLT's in Corollary  \ref{The Rate} hold true  at the rate $n^{-1}$.


\par
As to $G_{m_n}^{*}$ of (\ref{G^{*}}), on replacing $\big( \frac{w_{i}^{(n)}}{m_n}-\frac{1}{n}   \big)$  by $\big| \frac{w_{i}^{(n)}}{m_n}-\frac{1}{n} \big|$ in the proof of (a) of Corollary 2.1 of Cs\"{o}rg\H{o} \emph{et al.} \cite{Scand}, \emph{as $n,m_n \to +\infty$ so that} $m_n=o(n^2)$ we arrive at
\begin{equation}\label{(1.23)}
P_{X|w} (G_{m_n}^{*} \leq t) \to P(Z\leq t ) \ in \ probability-P_w\ for \ all \ t\in \mathds{R},
\end{equation}
and, via Lemma 1.2 in S. Cs\"{o}rg\H{o} and Rosalsky \cite{Csorgo and Rosalsky},   we conclude also the unconditional CLT
\begin{equation}\label{(1.24)}
P_{X,w} (G_{m_n}^{*} \leq t) \to P(Z\leq t ) \ for \ all \ t\in \mathds{R}.
\end{equation}
where, and also throughout, $Z$ stands for a standard normal random variable.

\par
For studying $G_{m_n}^{*}$  possibly in terms of  weights other than $w^{(n)}_{i}$, we refer to  Appendix 4.
\par

 \par

When $E_{X} X^2 <+\infty$, in   Appendix 1 we   show that when $n$ \emph{is fixed} and $m:=m_n \rightarrow +\infty$, the bootstrapped sample variance $S^{*^{2}}_{m_n}$,  as defined in (\ref{def. of S^*}),
converges in probability-$P_{X,w}$ to the sample variance $S_{n}^2$, i.e., \emph{as} $m:=m_n \rightarrow +\infty$ \emph{and} $n$ \emph{is fixed}, we have (cf. (\ref{(5.31)}) in Appendix 1 or Remark 2.1 of
Cs\"{o}rg\H{o} \emph{et al}. \cite{Scand})
\begin{equation}\label{retaining S^2 with fixed n}
S^{*^{2}}_{m} \rightarrow S_{n}^{2} \ in \ probability-P_{X,w}.
\end{equation}

\par
For related results along these lines in terms of $u$- and $v$-statistics,   we refer to Cs\"{o}rg\H{o} and Nasari \cite{J. multivariate Anal}, where in a more general setup, we establish in probability and almost sure consistencies  of bootstrapped $u$- and $v$-statistics.
\par
In   Appendix 1 we also show that, when $E_X X^2<+\infty$,  \emph{if} $n,\ m_n\rightarrow +\infty$ \emph{so that} $n=o(m_n)$,   then we have (cf. (\ref{(5.31)}) in Appendix 1)
\begin{equation}\label{bias of S*}
\big( S^{*^{2}}_{m_n} - S_{n}^{2}\big) \rightarrow 0 \ in \ probability-P_{X,w}.
\end{equation}
When $E_{X} X^4<+\infty$, the preceding convergence also holds true when $n=o(m^{2}_n)$ (cf. the proof of Corollary \ref{The Rate}).
 \par
 On combining (\ref{bias of S*}) with the CLT in (\ref{(1.24)}), when $E_X X^2<+\infty$, \emph{as} $n,m_n \to +\infty$ \emph{so that $m_n=o(n^2)$ and} $n=o(m_n)$, the following unconditional CLT holds true   as well in terms of $P_{X,w}$
\begin{equation}\label{(1.25)}
G_{m_n}^{**} \substack{d\\ \longrightarrow}\ Z ,
\end{equation}

where, and also throughout, ${\substack{d\\ \longrightarrow}}$ stands for convergence in distribution and  $G_{m_n}^{**}$ is as defined in (\ref{G^{**}}).

\par
Furthermore, when $E_X X^2 <+\infty$,  under the same conditions, i.e., on \emph{assuming that $n,m_n \to +\infty$ so that $m_n=o(n^2)$ and} $n=o(m_n)$, mutatis mutandis, via (b) of Lemma 2.1 of Cs\"{o}rg\H{o} \emph{et al.} \cite{Scand}, equivalently to (\ref{(1.25)}), in terms of $P_{X,w}$,  we also arrive at
\begin{equation}\label{(1.26)}
\tilde{G}_{m_n}^{**}:=\frac{\sum_{i=1}^n |\frac{w_{i}^{(n)}}{m_n}-\frac{1}{n}|(X_i-\mu)}{S_{m_n}^{*}/\sqrt{m_n} } \substack{d\\ \longrightarrow}\ Z .
\end{equation}
\par
Thus, without more ado, we may now conclude that the unconditional CLT's as in (\ref{(1.24)}), (\ref{(1.25)}) and (\ref{(1.26)}), respectively for $G_{m_n}^{*}$, ${G}_{m_n}^{**}$ and $\tilde{G}_{m_n}^{**}$, \emph{the bootstrapped versions of the traditional  Student pivot $T_{n}(X-\mu)$ for the population mean $\mu$} as in (\ref{equivalent to T_n}), can be used to construct exact size asymptotic C.I.'s  for the population mean $\mu=E_{X} X$ when $0< \sigma^2:=E_{X}(X-\mu)^2< +\infty$, \emph{if $n,m_n \to +\infty$ so that $m_n=o(n^2)$ in case of} (\ref{(1.24)}), and, \emph{in case of} (\ref{(1.25)}) \emph{and} (\ref{(1.26)}), \emph{so that $n=o(m_n)$ as well}. Moreover, when $E_{X} X^4 <+\infty$, \emph{if $n,m_n \to +\infty$ so that $m_n=o(n^2)$ }  and  $n=o(m^{2}_n)$ (cf. the proof of Corollary \ref{The Rate}), then    (\ref{(1.25)}) \emph{and} (\ref{(1.26)}) continue to hold true.
\par
We spell out the one based  on $G_{m_n}^{*}$ as in (\ref{(1.24)}). Thus, \emph{when $0< \sigma^2:=E_{X}(X-\mu)^2< +\infty$ and $m_n=o(n^2)$}, then,  for any $\alpha \in(0,1)$, we conclude (cf. also (\ref{(1.21)})-(\ref{doubleprime}) in comparison) a $1-\alpha$  size asymptotic C.I. for the population mean $\mu=E_{X} X$, which is also valid in terms of the conditional distribution $P_{X|w}$, as follows
\begin{equation}\label{(1.27)}
  \bar{X}_{n,m_n}^{*}-z_{\alpha/2}\frac{S_n \sqrt{.}}{\sum_{j=1}^n |\frac{w_{j}^{(n)}}{m_n}-\frac{1}{n}|}
\leq \mu \leq
 \bar{X}_{n,m_n}^{*}+z_{\alpha/2}\frac{S_n \sqrt{.}}{\sum_{j=1}^n |\frac{w_{j}^{(n)}}{m_n}-\frac{1}{n}|}
\end{equation}
where  $z_{\alpha/2}$ satisfies $P(Z \geq z_{\alpha/2})=\alpha/2$  and  $\sqrt{.}:= \sqrt{\sum_{j=1}^n (\frac{w_{j}^{(n)}}{m_n}-\frac{1}{n})^2}$.
\par
When $E_{X} X^2<+\infty$ on \emph{assuming} $n=o(m_n)$ and when $E_{X} X^4<+\infty$ on \emph{assuming} $n=o(m^{2}_n)$ (cf. the proof of Corollary \ref{The Rate}) as well, as $n,m_n \to +\infty$, then we can replace $S_n$ by $S_{m_n}^{*}$ in (\ref{(1.27)}), i.e., the thus obtained $1-\alpha$ size asymptotic C.I. for the population mean $\mu$ is then based on the CLT for ${G}_{m_n}^{**}$ as in (\ref{(1.25)}). A similar $1-\alpha$  size asymptotic C.I. for the population mean $\mu$ can be based on the CLT for $G^{**}_{m_n}$ as in (\ref{(1.26)}).

\section{\normalsize{Bootstrapped asymptotic pivots for the sample mean}}\label{section pivot for bar{X}}
We are now to present $T^{*}_{m_n}$, and further versions of it, as asymptotic bootstrapped   pivots for the sample mean $\bar{X_n}$ when $0<\sigma^2= E_{X}(X -\mu)^2<+\infty$. In view of its definition, when studying $T_{m_n}^{*}$, we may, \emph{without loss of generality}, forget about what the numerical value of  the population mean $\mu=E_X X$ could possibly be.
\par
To begin with, we consider the associated numerator term of $T_{m_n}^{*}$, and write
\\
\begin{eqnarray}
\sum_{i=1}^{n} (\frac{w^{(n)}_i}{m_n}-\frac{1}{n}) X_i &=& \frac{1}{m_n} \sum_{i=1}^{n} w^{(n)}_{i} X_i -\frac{1}{n}\sum_{i=1}^{n} X_i \label{bias}\\
&=& \bar{X^{*}}_{m_n}-\bar{X}_{n}. \nonumber
\end{eqnarray}
The term $\bar{X^{*}}_{m_n}$ is  known  as the bootstrapped sample  mean in the  literature.
\par
First we note that $E_{w|X}(\bar{X^{*}}_{m_n})=\bar{X}_{n}$, i.e., given $X$, the Efron bootstrap mean $\bar{X^{*}}_{m_n}$ is an unbiased estimator of the sample mean $\bar{X}_n$ and, consequently, (\ref{bias}) exhibits the bias of $\bar{X^{*}}_{m_n}$ \emph{vis-\`{a}-vis} $\bar{X}_{n}$, its conditional mean, given the data $X$. Moreover, when the original sample size $n$ \emph{is} assumed to be \emph{fixed}, then on taking \emph{only one} large bootstrap sub-sample of size $m:=m_n$, as $m \rightarrow +\infty$, we have
\begin{equation}\label{retaning bar{X} with fixed n}
\bar{X^{*}}_{m} \rightarrow  \bar{X}_n \ in \ probability-P_{X,w}
\end{equation}
(cf. (\ref{(5.30)}) of Appendix 1, where it is also shown that the bootstrap estimator of the sample mean that is based on taking a large number, $B$, of independent bootstrap sub-samples of size $n$ is equivalent to  taking \emph{only one}  large bootstrap sub-sample).
\par
Further to (\ref{retaning bar{X} with fixed n}), \emph{as} $n,\ m_n\rightarrow +\infty$, then (cf. (\ref{(5.30)})in Appendix 1)
\begin{equation}\label{bias of bar{X*}}
\big(\bar{X^{*}}_{m_n}- \bar{X}_n  \big) \rightarrow 0\ in\ probability-P_{X,w}.
\end{equation}

\par
Back to $T_{m_n}^{*}$ and further to (\ref{bias of bar{X*}}), we have that $E_{X|w}\big( \sum_{i=1}^n (\frac{w_{i}^{(n)}}{m_n}-\frac{1}{n}) X_i \big)=0$ and, \emph{if} $n,m_n \rightarrow +\infty$ \emph{so that} $m_n=o(n^2)$, then (cf. part (a) of Corollary 2.1 of Cs\"{o}rg\H{o} \emph{et al}. \cite{Scand})
\begin{equation}\label{(1.11)}
P_{X|w} (T_{m_n}^{*}\leq t) \rightarrow P(Z\leq t)\ in \ probability-P_w\ for \ all \ t\in \mathds{R},
\end{equation}
Consequently, as $n,m_n \to +\infty$ \emph{so that} $m_n=o(n^2)$, we arrive at (cf., e.g., Lemma 1.2 in S. Cs\"{o}rg\H{o} and Rosalsky \cite{Csorgo and Rosalsky})  
\begin{equation}\label{(1.12)}
P_{X,w} (T_{m_n}^{*}\leq t) \rightarrow P(Z\leq t)\ for \ all \ t\in \mathds{R},
\end{equation}
an unconditional CLT.
\par
Moreover, in view of the latter CLT and (\ref{bias of S*}), \emph{as} $n,m_n \rightarrow +\infty$ \emph{so that} $m_n=o(n^2)$ \emph{and}  $n=o(m_n)$, \emph{in terms of probability}-$P_{X,w}$ we conclude the unconditional CLT \begin{equation}\label{(1.13)}
{T}_{m_n}^{**} {\substack{d\\ \longrightarrow}}\ Z,
\end{equation}
where, $T_{m_n}^{**}$  is as defined in (\ref{T^{**}}). In turn we note also that, under the same conditions, i.e., \emph{on assuming that} $n,m_n \to +\infty$ \emph{so that} $m_n=o(n^2)$ \emph{and}  $n=o(m_n)$, via (b) of Corollary 2.1 of Cs\"{o}rg\H{o} \emph{et al} \cite{Scand}, equivalently to (\ref{(1.13)}) \emph{in terms of probability}-$P_{X,w}$ we also arrive at the unconditional CLT

\begin{equation}\label{(1.14)}
\tilde{T}_{m_n}^{**}:=  \frac{\sum_{i=1}^n (\frac{w_{i}^{(n)}}{m_n}-\frac{1}{n}) X_i}{S_{m_n}^{*}\big/ \sqrt{m_n}}  \ \substack{d\\ \longrightarrow}\ Z.
\end{equation}
\par
\par
Thus, when  \emph{conditioning on the weights}, via  (\ref{(1.11)}) we arrive at the  unconditional CLTs as in (\ref{(1.12)}), (\ref{(1.13)})  and (\ref{(1.14)}), respectively for $T_{m_n}^{*}$ and $\tilde{T}_{m_n}^{**}$  and $T_{m_n}^{**}$. These CLT's in hand for the  latter bootstrapped versions of the classical Student $t-$statistic $T_{n}(X)$ (cf. (\ref{def. of T_n}))  can be used to construct exact size asymptotic C.I.'s  for the sample mean $\bar{X}_n$ when $E_{X} X^2 <+\infty$, \emph{if $n,m_n \to +\infty$  so that $m_n=o(n^2)$ in case of} (\ref{(1.12)}) \emph{and, in case of} (\ref{(1.13)}) \emph{and} (\ref{(1.14)}), \emph{so that $n=o(m_n)$ as well}. Moreover, when $E_{X} X^4 <+\infty$, \emph{if $n,m_n \to +\infty$ so that $m_n=o(n^2)$ }  and  $n=o(m^{2}_n)$ (cf. the proof of Corollary \ref{The Rate}), then    (\ref{(1.13)}) \emph{and} (\ref{(1.14)}) continue to hold true.

\par
Consequently, under their respective conditions, any one of the unconditional CLT's in (\ref{(1.12)}), (\ref{(1.13)}) and (\ref{(1.14)})   can be used to construct exact size asymptotic confidence sets  for the sample mean $\bar{X}_n$,  \emph{treated as a random variable jointly with its bootstrapped  version} $\bar{X}_{m_n}^{*}$.
\par
We spell out the one based on $T_{m_n}^{*}$ as in (\ref{(1.12)}). Accordingly,  \emph{when $E_{X} X^2 <+\infty$ and $m_n,n \to +\infty$ so that $m_n=o(n^2)$}, then for any $\alpha \in (0,1)$, we conclude (cf. also  (\ref{bias})) a $1-\alpha$ size asymptotic confidence set, which is also valid in terms of the conditional distribution $P_{X|w}$,  that covers $\bar{X}_{n}$ as follows
\begin{equation}\label{(1.20)}
 \bar{X}_{m_n}^{*}-z_{\alpha/2} S_n \sqrt{.}\leq \bar{X}_n \leq \bar{X}_{m_n}^{*}+z_{\alpha/2} S_n \sqrt{.}  ,
\end{equation}
where $z_{\alpha/2}$ is as in (\ref{(1.27)}), and $\sqrt{.}:=\sqrt{\sum_{j=1}^n  (\frac{w_{i}^{(n)}}{m_n}- \frac{1}{n} )^2}$.

\par
When $E_{X} X^2<+\infty$, \emph{on assuming $n=o(m_n)$} and when $E_{X} X^4 <+\infty$ \emph{ on assuming $n=o(m^{2}_n)$  as well as} $m_n,n \to +\infty$, then we can replace $S_n$ by $S_{m_n}^{*}$ in (\ref{(1.20)}), i.e., then the thus obtained $1-\alpha$ asymptotic confidence set that covers $\bar{X}_n$ is based on the CLT for $\tilde{T}_{m_n}^{**}$  as in (\ref{(1.13)}). The equivalent CLT in (\ref{(1.14)}) can similarly be used to cover $\bar{X}_{n}$ when $E_{X} X^2 < +\infty$ and also when  $E_{X} X^4 <+\infty$.

\section{\normalsize{Bootstrapping a finite population when it is viewed  as a random sample of size $N$ from an infinite super-population}}\label{section super population}
As before,  let $X, X_1,\ldots$ be  independent random variables with a common distribution function $F$, mean $\mu=E_{X} X$ and, unless stated otherwise, variance $0< \sigma^2:=E_{X}(X-\mu)^2<+ \infty$. A super-population outlook regards a finite population as an imaginary large random sample of size $N$ from a hypothetical infinite population.  In our present context,  we consider $\{X_1, \ldots, X_N \}$ to be a \emph{concrete} or \emph{imaginary random sample} of size $N\geq 1$ on $X$, and study it as a finite population of $N$ independent and identically distributed random variables with respective mean and variance
\begin{equation}\label{(1.28)}
\bar{X}_{N}:= \sum_{i=1}^N X_i /N \ \textrm{and} \ S_{N}^2:=\sum_{i=1}^N (X_i-\bar{X}_N)^{2}/N
\end{equation}
\emph{that are to be estimated} via taking samples from it.
\par
Our approach in this regard is to assume that we can view an  \emph{imaginary random sample} $\{X_1, \ldots, X_N \}$ as a \emph{finite population of real valued random variables  with $N$ labeled units},  and sample its set of indices $\{1, \ldots, N \}$ with replacement $m_N$ times so that for each $1\leq i \leq N$, $w_{i}^{(N)}$ is the count of the number of times the index $i$ of $X_i$ is chosen in this re-sampling process. Consequently, as to the weights $w_{i}^{(N)}$, $1\leq i \leq N$, mutatis mutandis, Remark \ref{Remark 1} obtains with $m_N=\sum_{i=1}^N w_{i}^{(N)}$ and, for each $N\geq 1$, the multinomial weights
$$(w^{(N)}_{1},\ldots,w^{(N)}_{N})\  \substack{d\\=}\ \ multinomial(m_{N};\frac{1}{N},\ldots,\frac{1}{N})
$$
are independent from the finite population of the $N$ labeled units $\{X_1, \ldots, X_N \}$.
\par
In the just described context, mutatis mutandis, under their respective conditions the respective  results  as in (\ref{retaining S^2 with fixed n}) and (\ref{retaning bar{X} with fixed n})   are applicable to consistently estimate $\bar{X}_{N}$ and $S_{N}^{2}$ of (\ref{(1.28)}) respectively with $N$ fixed and $m:=m_N \to +\infty$, while  the respective CLT's in (\ref{(1.13)}) and  (\ref{(1.14)}) can be adapted to construct exact size asymptotic confidence sets for covering $\bar{X}_{N}$. For illustrating this, \emph{\`{a} la} ${T}_{m_n}^{**}$ in (\ref{(1.13)}), we let
\\
\begin{eqnarray}
{T}_{m_N}^{**}&:=& \frac{\sum_{i=1}^N (\frac{w_{i}^{(N)}}{m_N}-\frac{1}{N}) X_i}{S_{m_N}^{*} \sqrt{\sum_{j=1}^N (\frac{w_{j}^{(N)}}{m_N}-\frac{1}{N})^2} } \nonumber\\
&=:& \frac{\bar{X}_{m_N}^{*}-\bar{X_N}}{S_{m_N}^{*} \sqrt{\sum_{j=1}^N (\frac{w_{j}^{(N)}}{m_N}-\frac{1}{N})^2} },\label{(1.30)}
\end{eqnarray}
where $\bar{X}_{m_N}^{*}:=\sum_{i=1}^N  X_i w_{i}^{(N)}/m_N$, the bootstrapped finite population mean   and, \emph{\`{a} la} $S_{m_n}^{*}$ in (\ref{def. of S^*}),
\begin{equation}\label{(1.31)}
S^{*^{2}}_{m_N}:= \sum_{i=1}^N w_{i}^{(N)} \big( X_i- \bar{X^{*}}_{m_{N}} \big)^2 \big/ m_N
\end{equation}
is the bootstrapped finite population variance.
\par
With $N$ \emph{fixed and} $m=m_N\to +\infty$, when $X$ has a finite variance, \emph{\`{a} la} (\ref{retaning bar{X} with fixed n}) and (\ref{retaining S^2 with fixed n}) respectively, we arrive at
\begin{equation}\label{(1.32)}
\bar{X^{*}}_{m} \rightarrow  \bar{X}_N \ in \ probability-P_{X,w}
\end{equation}
and
\begin{equation}\label{(1.33)}
S^{*^{2}}_{m} \rightarrow S_{N}^{2} \ in \ probability-P_{X,w},
\end{equation}
i.e.,  $\bar{X^{*}}_{m}$ and $S^{*^{2}}_{m}$ are consistent estimators of the finite population mean $\bar{X_N}$ and variance $S_{N}^{2}$ respectively. Also, \emph{\`{a} la} (\ref{bias of bar{X*}}), \emph{as} $N, m_N \to +\infty$, we conclude
\begin{equation}\label{(1.34)}
\big(\bar{X^{*}}_{m_N}- \bar{X}_N  \big) \rightarrow 0\ in\ probability-P_{X,w}
\end{equation}
and, \emph{if} $N,\ m_N\rightarrow +\infty$ \emph{so that} $N=o(m_N)$,   then \emph{\`{a} la} (\ref{bias of S*}), we also have
\begin{equation}\label{(1.35)}
\big( S^{*^{2}}_{m_N} - S_{N}^{2}\big) \rightarrow 0 \ in \ probability-P_{X,w}.
\end{equation}
\par
Furthermore, when studying ${T}_{m_N}^{**}$ as in (\ref{(1.30)})  and, \emph{\`{a} la} $\tilde{T}_{m_n}^{**}$ as in (\ref{(1.14)}),
\begin{equation}\label{(1.36)}
\tilde{T}_{m_N}^{**}:=  \frac{\sum_{i=1}^N (\frac{w_{i}^{(N)}}{m_N}-\frac{1}{N}) X_i}{S_{m_N}^{*}\big/ \sqrt{m_N}},
\end{equation}
we may,  without loss of generality, forget about what the value of  the super-population mean $\mu=E_X  X$ would be like.   On assuming that $0<E_{X} X^2<+\infty$, \emph{if} $N,m_N \to +\infty$ \emph{so that} $m_N=o(N^2)$ \emph{and} $N=o(m_N)$, then in view of the respective equivalent statements of (\ref{(1.13)}) and (\ref{(1.14)}), \emph{in terms of probability-$P_{X,w}$}, we arrive at having the unconditional CLT's

\begin{equation}\label{(1.38)}
T_{m_N}^{**} {\substack{d\\ \longrightarrow}}\ Z
\end{equation}

and

\begin{equation}\label{(1.37)}
\tilde{T}_{m_N}^{**} {\substack{d\\ \longrightarrow}}\ Z.
\end{equation}

\par
Back to ${T}_{m_N}^{**}$ as in (\ref{(1.30)}), via the CLT of (\ref{(1.37)}), for any $\alpha\in(0,1)$, we conclude a $1-\alpha$ size asymptotic confidence set, which is also valid in terms of the conditional distribution $P_{X|w}$,   that covers $\bar{X}_N$ as follows
\begin{equation}\label{(1.39)}
 \bar{X}_{m_N}^{*}-z_{\alpha/2}S_{m_N}^{*} \sqrt{.}
\leq \bar{X}_N \leq
 \bar{X}_{m_N}^{*}+z_{\alpha/2}S_{m_N}^{*} \sqrt{.}
\end{equation}
with $z_{\alpha/2}$ as in (\ref{(1.27)}) and $\sqrt{.}:= \sqrt{\sum_{j=1}^N (\frac{w_{j}^{(N)}}{m_N}-\frac{1}{N})^2}$.

\par
In the present context of taking samples with replacement from a finite population with $N$ labeled units, $\{X_1,\ldots,X_N \}$,  that is viewed as an imaginary random sample of size $N$ from an infinite super-population with mean $\mu=E_{X} X$ and variance $0<\sigma^2=E_{X}(X-\mu)^2<+\infty$, we can also estimate the population mean $\mu$ via adapting the results of of our Section \ref{section pivot for mu} to fit this setting.
\par
To begin with, the bootstrapped weighted average (cf. (\ref{(1.21)})) in this context becomes
\begin{equation}\label{(1.40)}
\bar{X}_{m_N}^{*}(\mu):= {\sum_{i=1}^N |\frac{w_{i}^{(N)}}{m_N}-\frac{1}{N}| (X_i-\mu)}
\end{equation}
and, given the weights, 
we have
\begin{equation}\label{(1.41)}
E_{X|w} (\bar{X}_{m_N}^{*}(\mu))=0.
\end{equation}
\par
Furthermore, with the finite population $N$ \emph{fixed} and  $m:=m_N$, \emph{as} $m \to +\infty$, we also conclude (cf. (\ref{(1.22)}))
\begin{equation}\label{(1.42)}
\bar{X}_{m_{N}}^{*}(\mu)=\bar{X}_{m}^{*}(\mu) \to 0 \ in \ probability-P_{X,w},
\end{equation}
and the same holds true if $N \to +\infty$ as well.
\par
A consistent estimator of the super population mean $\mu$, similarly to (\ref{(1.22)'}), can be defined as
\begin{equation}\label{(1.22)'N}
\frac{\sum_{i=1}^N |\frac{w_{i}^{(N)}}{m_N}-\frac{1}{N}| X_i}{\sum_{j=1}^N |\frac{w_{j}^{(N)}}{m_N}-\frac{1}{N}|}=: \bar{X}_{N,m_N}^{*}.
\end{equation}
Thus, similarly to (\ref{doubleprime}), it can be shown that, as  $N,m_N \to +\infty$ such that $m_N=o(N^2)$, $\bar{X}_{N,m_N}^{*}$ is a consistent estimator for the super-population mean $\mu$, in terms of $P_{X,w}$, i.e.,
\begin{equation}\label{(1.22)''}
\bar{X}_{N,m_N}^{*} \to \mu \ in \ probability-P_{X,w}.
\end{equation}

\par
Now, with the bootstrapped finite population variance  $S_{m_N}^{*^{2}}$ as in (\ref{(1.31)}) (cf. also (\ref{(1.35)})), along the lines of arguing the conclusion of (\ref{(1.25)}), \emph{as} $N,m_N \to +\infty$ \emph{so that} $m_N=o(N^2)$ \emph{and} $N=o(m_N)$, we conclude the following unconditional  CLT as well
\begin{equation}\label{(1.43)}
\tilde{G}_{m_N}^{**}:=\frac{\sum_{i=1}^N |\frac{w_{i}^{(N)}}{m_N}-\frac{1}{N}|(X_i-\mu)}{S_{m_N}^{*} \sqrt{\sum_{j=1}^N (\frac{w_{j}^{(N)}}{m_N}-\frac{1}{N})^2}} \substack{d\\ \longrightarrow}\ Z .
\end{equation}
Furthermore, under the same conditions, equivalently to (\ref{(1.43)}), \emph{\`{a} la} (\ref{(1.26)}), we also have the unconditional CLT  in parallel to that of (\ref{(1.43)})
\begin{equation}\label{(1.44)}
G_{m_N}^{**}:=\frac{\sum_{i=1}^N |\frac{w_{i}^{(N)}}{m_N}-\frac{1}{N}|(X_i-\mu)}{S_{m_N}^{*}/\sqrt{m_N} } \substack{d\\ \longrightarrow}\ Z .
\end{equation}
\par
Both (\ref{(1.43)}) and (\ref{(1.44)}) can be used to construct a $1-\alpha$ size asymptotic C.I., with $\alpha\in (0,1)$, for the \emph{unknown} super-population mean $\mu$. We spell out the one based on the CLT of (\ref{(1.43)}). Accordingly, \emph{as $N, m_N \to +\infty$ so that $m_N=o(N^2)$ and $N=o(m_N)$}, we arrive at the following $1-\alpha$ size asymptotic C.I., which is valid in terms of the conditional distribution $P_{X|w}$,  as well as in terms of the   joint distribution $P_{X,w}$,  for the \emph{unknown} super-population mean $\mu$
\begin{equation}\label{(1.45)}
 \bar{X}_{N,m_N}^{*}-z_{\alpha/2}\frac{S^{*}_{m_N} \sqrt{.}}{\sum_{j=1}^N |\frac{w_{j}^{(N)}}{m_N}-\frac{1}{N}|}
\leq \mu \leq
 \bar{X}_{N,m_N}^{*}+z_{\alpha/2}\frac{S^{*}_{m_N} \sqrt{.}}{\sum_{j=1}^N |\frac{w_{j}^{(N)}}{m_N}-\frac{1}{N}|}
\end{equation}
with $z_{\alpha/2}$ as in (\ref{(1.27)}) and $\sqrt{.}:= \sqrt{\sum_{j=1}^N (\frac{w_{j}^{(N)}}{m_N}-\frac{1}{N})^2}$, a companion of the $1-\alpha$ size asymptotic confidence set of (\ref{(1.39)}) that, in the same super-population setting, covers the \emph{unknown} finite population mean $\bar{X}_N$ under the same conditions.

\section{\normalsize{Bootstrapped CLT's and C.I.'s for the empirical and theoretical distributions}}\label{CLT empirical}
Let $X,X_1,X_2, \ldots$ be independent real valued random variables with a common distribution function $F$ as before, but without assuming the existence of any finite moments for $X$. Consider $\{ X_1,\ldots,X_N\}$ a \emph{concrete  or imaginary random sample} of size $N\geq 1$ on $X$ of a hypothetical infinite population, and defined their empirical distribution function
\begin{equation}\label{(5.1)}
F_{N}(x):= \sum_{i=1}^N \mathds{1}(X_i \leq x)/N, \ x\in \mathbb{R},
\end{equation}
and the sample variance of the indicator variable $\mathds{1}(X_i \leq x)$ by
\begin{equation}\label{(5.2)}
S^{2}_{N}:=\frac{1}{N} \sum_{i=1}^N \big( \mathds{1}(X_i \leq x) - F_{N}(x)\big)^2= F_N (x)(1-F_N (x)),\ x\in\mathbb{R},
\end{equation}
\emph{that, together with the theoretical distribution function $F$, are to be estimated} via taking samples from $\{X_1,\ldots, X_N \}$  as in our previous sections in general, and as in
Section \ref{section super population} in particular. On replacing $N$ by $n$ in case of a real sample of size $n$, our present general general formulation of the problems in hand, as well as the results thus concluded, mutatis mutandis, continue to hold true when interpreted in the context of Sections \ref{section pivot for mu} and  \ref{section pivot for bar{X}}   that deal with a concrete  random sample $\{X_1,\dots,X_n \}$ of size $n\geq 1$.
\par
Accordingly, we view a  \emph{concrete} or \emph{imaginary random sample } $\{X_1,\ldots,X_N\}$ as a \emph{finite population} of \emph{real valued random variables} with $N$ label units, and sample its set of indices $\{1,\ldots,N \}$ with replacement $m_N$ times so that for $1\leq i \leq N$, $w^{(N)}_{i}$ is the count of the  number of times the index $i$ of $X_i$ is chosen in this re-sampling process and, as in Section \ref{section super population}, $m_N=\sum_{i=1}^{N} w^{(N)}_{i} $ and, for each $N \geq 1$,  the multinomial weights
\begin{equation*}
\big(w^{(N)}_{1},\ldots,w^{(N)}_{N}\big) {\substack{d\\=}} \ multinomial \big(m_N;\frac{1}{N},\ldots,\frac{1}{N}\big)
\end{equation*}
are independent from the finite population of the $N$ labeled units $\{X_1,\ldots,X_N \}$.

Define the standardized bootstrapped  empirical process
\begin{eqnarray}
\alpha^{(1)}_{m_{N},N}(x)&:=& \frac{\sum_{i=1}^N (\frac{w^{(N)}_{i}}{m_N}-\frac{1}{N}  ) \mathds{1}(X_i \leq x)}{\sqrt{F(x)(1-F(x))} \sqrt{\sum_{j=1}^N (\frac{w^{(N)}_{j}}{m_N}-\frac{1}{N}  )^2} } \label{(5.3)}\\
&=& \frac{\sum_{i=1}^N \frac{w^{(N)}_{i}}{m_N} \mathds{1}(X_i \leq x)-F_{N}(x)}{\sqrt{F(x)(1-F(x))} \sqrt{\sum_{j=1}^N (\frac{w^{(N)}_{j}}{m_N}-\frac{1}{N}  )^2}}\nonumber\\
&=& \frac{F^{*}_{m_{N},N}(x) -F_{N}(x)  }{\sqrt{F(x)(1-F(x))} \sqrt{\sum_{j=1}^N (\frac{w^{(N)}_{j}}{m_N}-\frac{1}{N}  )^2}}, \ x\in \mathbb{R}\nonumber
\end{eqnarray}
where
\begin{equation}\label{(5.3)}
F^{*}_{m_{N},N}(x):= \sum_{i=1}^N \frac{w^{(N)}_{i}}{m_N} \mathds{1}(X_i \leq x), \ x\in \mathbb{R},
\end{equation}
is the bootstrapped empirical process.
\par
We note that
\begin{equation}\label{(5.4)}
E_{X|w}(F^{*}_{m_{N},N}(x))=F(x), \ E_{w|X}(F^{*}_{m_{N},N}(x))=F_N (x)\ and \ E_{X,w}(F^{*}_{m_{N},N}(x))=F(x).
\end{equation}
 \par
 Define also the bootstrapped finite population variance of the indicator random variable $\mathds{1}(X_i \leq x)$ by putting
\begin{eqnarray}
S^{*^{2}}_{m_{N},N}(x)&:=& \sum_{i=1}^{N} w^{(N)}_{i} \big( \mathds{1}(X_i \leq x)- F^{*}_{m_{N},N}(x)\big)^2\big/m_N \label{(5.4)}\\
&=& F^{*}_{m_{N},N}(x)(1-F^{*}_{m_{N},N}(x)), \ x\in \mathbb{R}.\nonumber
\end{eqnarray}

\par
With $N$ \emph{fixed and } $m=m_N\rightarrow +\infty$, along the lines of (\ref{(1.32)}) we arrive at

\begin{equation}\label{(5.6)}
F^{*}_{m_{N},N}(x)\longrightarrow F_{N}(x)\ in\ probbility-P_{X,w}\ pointwise \ in \ x\in \mathbb{R},
\end{equation}
and, consequently, point-wise in $x \in \mathbb{R}$, \emph{as} $m=m_N\rightarrow +\infty$,
\begin{equation}\label{(5.7)}
S^{*^{2}}_{m_{N},N}(x) \longrightarrow F_N(x)(1-F_N (x))=S^{2}_{N}(x)\ in \ probability-P_{X,w}.
\end{equation}
Furthermore, \emph{\`{a} la} (\ref{(1.34)}), \emph{as} $N,m_N\to +\infty$, \emph{pointwise in} $x \in \mathbb{R}$, we conclude
\begin{equation}\label{(5.8)}
\big(F^{*}_{m_{N},N}(x)-F_N (x)\big)\longrightarrow 0\ in\ probbility-P_{X,w,}
\end{equation}
that, in turn, \emph{poitwise in} $x\in \mathbb{R}$, \emph{as}  $N,m_N\to +\infty$, implies
\begin{equation}\label{(5.9)}
\big(S^{*^{2}}_{m_{N},N}-S^{2}_{N}(x) \big) \longrightarrow 0=\ in \ probability-P_{X,w}.
\end{equation}

\par
We wish to note and emphasize that, unlike in (\ref{(1.35)}), for concluding (\ref{(5.9)}), \emph{we do not assume that} $N=o(m_N)$ as $N,m_N\to +\infty$.
\par

Further to the standardized bootstrapped empirical process $\alpha^{(1)}_{N,m_N}(x)$, we now define the following Studentized/self-normalized bootstrapped versions of this process:

\begin{eqnarray}
\hat{\alpha}^{(1)}_{m_{N},N}(x)&:=& \frac{\sum_{i=1}^N (\frac{w^{(N)}_{i}}{m_N}-\frac{1}{N}  ) \mathds{1}(X_i \leq x)}{\sqrt{F_N(x)(1-F_N(x))} \sqrt{\sum_{j=1}^N (\frac{w^{(N)}_{j}}{m_N}-\frac{1}{N}  )^2} } \label{(5.10)}\\
\hat{\hat{\alpha}}^{(1)}_{m_{N},N}(x)&:=& \frac{\sum_{i=1}^N (\frac{w^{(N)}_{i}}{m_N}-\frac{1}{N}  ) \mathds{1}(X_i \leq x)}{\sqrt{F^{*}_{m_{N},N}(x)(1-F^{*}_{m_{N},N}(x))} \sqrt{\sum_{j=1}^N (\frac{w^{(N)}_{j}}{m_N}-\frac{1}{N}  )^2} } \label{(5.11)} \\
\hat{\alpha}^{(2)}_{m_{N},N}(x)&:=& \frac{\sum_{i=1}^N \big|\frac{w^{(N)}_{i}}{m_N}-\frac{1}{N} \big| \big(\mathds{1}
(X_i \leq x)-F(x)\big)}{\sqrt{F_N(x)(1-F_N(x))} \sqrt{\sum_{j=1}^N (\frac{w^{(N)}_{j}}{m_N}-\frac{1}{N}  )^2} } \label{(5.12)} \\
\hat{\hat{\alpha}}^{(2)}_{m_{N},N}(x)&:=& \frac{\sum_{i=1}^N \big|\frac{w^{(N)}_{i}}{m_N}-\frac{1}{N}  \big| \big(\mathds{1}(X_i \leq x)-F(x)\big)}{\sqrt{F^{*}_{m_{N},N}(x)(1-F^{*}_{m_{N},N}(x))} \sqrt{\sum_{j=1}^N (\frac{w^{(N)}_{j}}{m_N}-\frac{1}{N}  )^2} }.\label{(5.13)}
\end{eqnarray}
\par
Clearly, on replacing $X_i$ by $\mathds{1}(X_i \leq x)$ and $\mu$ by $F(x)$ in the formula  in  (\ref{(1.40)}), we arrive at the respective  statements of (\ref{(1.41)}) and (\ref{(1.42)}) in this context. Also, replacing  $X_i$ by $\mathds{1}(X_i \leq x)$ in the formula as in (\ref{(1.22)'N}), we conclude the statement of (\ref{(1.22)''}) with $\mu$ replaced $F(x)$.
\par
In Lemma 5.2 of Cs\"{o}rg\H{o} \emph{et al}. \cite{Scand} it is shown that, \emph{if} $m_N , N\to +\infty$ \emph{so that $m_N=o(N^2)$, then}

\begin{equation}\label{(5.14)}
M_N:=\frac{ \max_{1\leq i \leq N} \big( \frac{w^{(N)}_{i}}{N}-\frac{1}{N}\big)^2 } {\sum_{j=1}^N \big( \frac{w^{(N)}_j}{N}-\frac{1}{N}\big)^2 }\to 0\ in \ probability-P_{w}.
\end{equation}
This, mutatis mutandis, combined with (a) of Corollary 2.1 of Cs\"{o}rg\H{o} \emph{et al}. \cite{Scand}, \emph{as $N, m_N \to +\infty$ so that} $m_N=o(N^2)$, yields
\begin{equation}\label{(5.15)}
P_{X|w} \big( \hat{\alpha}^{(s)}_{m_{N},N}(x)\leq t \big)\to P(Z \leq t) \ in \ probability-P_{w},\ for\ all \ x,t\in \mathbb{R},
\end{equation}
with $s=1$ and also for $s=2$, and via  Lemma 1.2 in  S. Cs\"{o}rg\H{o} and Rosalsky \cite{Csorgo and Rosalsky}, this results in having also the unconditional CLT

\begin{equation}\label{(1.16)}
P_{X,w} \big( \hat{\alpha}^{(s)}_{m_{N},N}(x)\leq t \big)\to P(Z \leq t) \ for\ all \ x,t\in \mathbb{R},
\end{equation}
with $s=1$ and also for $s=2$.

\par
On combining (\ref{(1.16)}) and (\ref{(5.9)}), \emph{as } $N,m_N\to +\infty$ \emph{so that} $m_N=o(N^2)$, when $s=1$ in (\ref{(1.16)}),  we  conclude

\begin{equation}\label{(5.17)}
\hat{\hat{\alpha}}^{(1)}_{m_{N},N}(x)
\substack{d\\ \longrightarrow} Z
\end{equation}
and, when $s=2$ in (\ref{(1.16)}), we arrive at
\begin{equation}\label{(5.18)}
\hat{\hat{\alpha}}^{(2)}_{m_{N},N}(x)
\substack{d\\ \longrightarrow} Z
 \end{equation}
for all $x \in \mathbb{R}$.

\par
Furthermore, \emph{as}  $m_N, N \to +\infty$ \emph{so that} $m_N=o(N^2)$,  in addition to (\ref{(5.14)}), we also have (cf. (\ref{for appendix 2}))

\begin{equation}\label{(5.19)}
\big| \sum_{i=1}^{N} \big(\frac{w^{(N)}_i }{N} -\frac{1}{N} \big)^2-\frac{1}{m_N}  \big|=o_{P_{w}}(1).
\end{equation}
Consequently, in the CLT's of (\ref{(5.15)})-(\ref{(5.18)}), the term $\sqrt{\sum_{i=1}^{N} \big(\frac{w^{(N)}_i }{N} -\frac{1}{N} \big)^2}$ in the respective denumerators of $\hat{\alpha}^{(s)}_{m_{N},N}(x)$ and $\hat{\hat{\alpha}}^{(s)}_{m_{N},N}(x) $, $s=1,2$, can be replaced by $1/\sqrt{m_N}$.

\par
In case of a \emph{concerted random samples} of size $N\geq 1$ on $X$, the CLT's for $\hat{\alpha}^{(1)}_{m_{N},N}(x) $ and $\hat{\hat{\alpha}}^{(1)}_{m_{N},N}(x)$ can both be used to construct a $1-\alpha$,  $\alpha \in (0,1)$, size asymptotic confidence sets for covering the empirical distribution function $F_{N} (x)$, point-wise in $x \in \mathbb{R}$, as $N,m_N \to +\infty$ \emph{so that} $m_N =o(N^2)$, while in case of an \emph{imaginary random sample}  of size $N\geq 1$ on $X$ of a hypothetical infinite population, the CLT for $\hat{\hat{\alpha}}^{(1)}_{m_{N},N}(x)$ works also similarly estimating the, in this case, unknown empirical distribution function $F_{N} (x)$. The respective CLT's for $\hat{\alpha}^{(2)}_{m_{N},N}(x)$  and $\hat{\hat{\alpha}}^{(2)}_{m_{N},N}(x)$ work in a similar way for point-wise estimating the \emph{unknown} distribution function $F(x)$ of a \emph{hypothetical infinite population} in case of a \emph{concrete} and \emph{imaginary random sample} of size $N \geq 1$ on $X$. Furthermore, mutatis mutandis, the Berry-Esse\'{e}n type inequality  (A) of our Theorem \ref{Berry-Esseen} continues to hold true for both $\hat{\alpha}^{(2)}_{m_{N},N}(x)$ and $\hat{\hat{\alpha}}^{(2)}_{m_{N},N}(x)$, and  so does  also (B)   of Theorem \ref{Berry-Esseen} for both  $\hat{\alpha}^{(1)}_{m_{N},N}(x)$  and $\hat{\hat{\alpha}}^{(1)}_{m_{N},N}(x)$, without the assumption $E_{X} |X|^{3}<+\infty$, for the indicator random variable $\mathds{1}(X \leq x)$ requires no moments assumptions.

\begin{remark}
In the context of this section, (A) and (B) of Corollary \ref{The Rate} read as follows: As $N,m_N \to +\infty$ in such a way that $m_N=o(N^2)$, then, mutatis mutandis, (A) and (B) hold true for  $\hat{\alpha}^{(1)}_{m_{N},N}(x)$ and $\hat{\hat{\alpha}}^{(1)}_{m_{N},N}(x)$, with $O(\max\{m_N/N^{2}, 1/m_N\})$
in both.  Consequently, on taking $M_N=N$, we immediately obtain the optimal $O(N^{-1})$ rate conclusion of Remark \ref{Remark 2.1} in this context, i.e., uniformly in $t \in \mathbb{R}$ and point-wise in $x \in \mathbb{R}$ for $\hat{\alpha}^{(1)}_{m_{N},N}(x)$ and ${\hat{\alpha}}^{(2)}_{m_{N},N}(x)$.

Also, for a rate of convergence  of the respective CLT's  via (C) and  (D) of Corollary \ref{The Rate} for  $\hat{\hat{\alpha}}^{(1)}_{m_{N},N}(x)$ and  $\hat{\hat{\alpha}}^{(2)}_{m_{N},N}(x)$,  as $N,m_N \to +\infty$ in such away that $m_N=O(N^2)$,    we obtain the rate $O(\max\{m_N/N^{2}, 1/m_N\})$. The latter means that observing a sub-sample of   size $m_N=N^{1/2}$ will result in the rate of convergence of $N^{-1/2}$ for $\hat{\hat{\alpha}}^{(1)}_{m_{N},N}(x)$ and  $\hat{\hat{\alpha}}^{(2)}_{m_{N},N}(x)$,  uniformly in $t \in \mathbb{R}$ and point-wise in $x \in \mathbb{R}$.

\end{remark}

\par
In the context of this section, asymptotic   $1-\alpha$ size confidence sets, in terms of the joint distribution $P_{X,w}$ which are also valid in terms of $P_{X|w}$,  for the empirical and the population distribution functions $F_{N}(x)$ and  $F(x)$,  for each $x \in \mathbb{R}$, are in fact of the forms (\ref{(1.39)}) and (\ref{(1.45)}), respectively,  on replacing $X_i$ by $\mathds{1}(X_i\leq x)$ which are spelled out as follows

\begin{equation*}
 F^{*}_{m_{N},N}(x)-z_{\alpha/2}S_{m_N}^{*} \sqrt{.}
\leq F_{N}(x) \leq
 F^{*}_{m_{N},N}(x)+z_{\alpha/2}S_{m_N}^{*} \sqrt{.}
\end{equation*}

\begin{equation*}
 F^{*}_{m_{N},N}(x) - z_{\alpha/2}\frac{S^{*}_{m_N} \sqrt{.}}{\sum_{j=1}^N |\frac{w_{j}^{(N)}}{m_N}-\frac{1}{N}|}
\leq F(x) \leq
 F^{*}_{m_{N},N}(x) + z_{\alpha/2}\frac{S^{*}_{m_N} \sqrt{.}}{\sum_{j=1}^N |\frac{w_{j}^{(N)}}{m_N}-\frac{1}{N}|}
\end{equation*}
with $z_{\alpha/2}$ as in (\ref{(1.27)}),  $\sqrt{.}:= \sqrt{\sum_{j=1}^N (\frac{w_{j}^{(N)}}{m_N}-\frac{1}{N})^2}$ and $S_{m_N}^{*}=F^{*}_{m_{N},N}(x)(1-F^{*}_{m_{N},N}(x))$.

\section{\normalsize{Comparison to classical bootstrap C.I.'s}}\label{Comparison section}
In this section we show that, unlike our direct pivot $G^{*}_{m_n}$ of (\ref{G^{*}}) for $\mu$, the indirect use of  $T^{*}_{m_n}$ of (\ref{T^{*}}) along the lines of  the classical method of constructing a bootstrap C.I.  for the population mean $\mu$ does not lead to  a better error rate  than that of the classical CLT for $T_{n}(X-\mu)$, which is $n^{-1/2}$ under the conditions of Corollary \ref{The Rate}.
\par
A bootstrap estimator for a quantile $\eta_{n,1-\alpha}$, $0< \alpha <1$, of $T_{n}(X-\mu)$ is the solution of the inequality
\begin{equation}\label{(6.1)}
P(T^{*}_{m_n}\leq x| X_1,\ldots,X_n):=P_{w|X}(T^{*}_{m_n}\leq x )\geq 1-\alpha,
\end{equation}
i.e., the smallest value of $x=\hat{\eta}_{n,1-\alpha}$ that satisfies (\ref{(6.1)}), where $T^{*}_{m_{n}}$ is as defined in (\ref{T^{*}}).  Since the latter bootstrap quantiles should be close to the true quantiles $\eta_{n,1-\alpha}$ of $T_{n}(X-\mu)$, in view of (\ref{(6.1)}), it should be true that
\begin{equation}\label{(6.2)}
P_{w|X}(T_n (X-\mu)\leq \hat{\eta}_{n,1-\alpha} )\approx 1-\alpha.
\end{equation}

\par
In practice, the value of $\hat{\eta}_{n,\alpha}$ is usually estimated by simulation (cf. Efron and Tibshirani \cite{Efron and Tibshirani}, Hall \cite{Hall B}) via producing $B\geq 2$ independent copies of $T^{*}_{m_n}$,  usually  with $m_n=n$,  given $X_1,\ldots,X_n$.

\par
A classical bootstrap C.I. of level $1-\alpha$, $0 < \alpha <1$, for $\mu$ is constructed by using $T_{n} (X-\mu)$ as a pivot and estimating the cutoff point $\hat{\eta}_{n,1-\alpha}$ using the  $(B+1).(1-\alpha)$th largest value  of $T^{*}_{m_n}(b)$, $1\leq b \leq B$, where,   each $T^{*}_{m_n}(b)$, is computed based on the $b$-th bootstrap sub-sample.
  We note that the preceding method of constructing a classical bootstrap C.I. at level  $1-\alpha$ is for the case when  $(B+1).(1-\alpha)$ is an integer already.
\par
For the sake of comparison of  our main results,  namely Theorem \ref{Berry-Esseen} and Corollary \ref{The Rate},  to the classical bootstrap C.I.'s which were also investigated by
 Hall  \cite{Hall B}, here we study  the rate of convergence of classical  bootstrap C.I.'s.
 \par
 As it will be seen below, our investigations agree with those of Hall \cite{Hall B} in   concluding  that the number of bootstrap sub-samples $B$ does not have to be particularly large for a classical  bootstrap C.I. to reach its nominal probability  coverage.
   \par
   In this section,  we also show  that the rate at which  the  probability of the event that  the conditional probability, given $w^{(n)}_i$'s, of  a classical bootstrap C.I. for $\mu$  deviating from its nominal probability coverage  by  any given positive number, vanishes  at a   rate  that can, at best, be $O(n^{-1/2})$, as $n\to +\infty$ (cf. Theorem \ref{rate of bootstrap}). This is in contrast to  our  Corollary \ref{The Rate} and Remark \ref{Remark 2.1}. It is also noteworthy  that the rate of convergence of a joint distribution is essentially the same as
    those  of its conditional versions to the same limiting distribution. Therefore, the  preceding rate of, at best, $O(n^{-1/2})$ for the conditional, given $w^{(n)}_i$'s,  probability coverage of  classical bootstrap C.I.'s is  inherited by their probability coverage in terms of the joint distribution.

 \par
In what follows    $B\geq 2$ is a \emph{fixed} positive integer, that stands  for  the number of bootstrap sub-samples of size $m_n$, generated via $B$ times independently re-sampling with replacement from $\{X_i, 1\leq i \leq n \}$. We  let $T^{*}_{m_n}(1),\ldots,T^{*}_{m_n}(B)$ be the   versions of computed $T^{*}_{m_n}$ based on the drawn $B$ bootstrap sub-samples. We now   state  a multivariate CLT  for the $B$ dimensional vector
\begin{equation}\nonumber
\big( \frac{T_{n}(X-\mu)-T^{*}_{m_n}(1)}{\sqrt{2}},\ldots, \frac{T_{n}(X-\mu)-T^{*}_{m_n}(B)}{\sqrt{2}} \big)^T.
\end{equation}
 The reason for  investigating the asymptotic distribution of the latter vector has to do with computing the actual probability  coverage of the classical bootstrap C.I.'s as in (\ref{bootstrap CI}) below. Furthermore, in Section \ref{Section 6.2}, we compute
the actual probability coverage of the classical  bootstrap C.I.'s as in (\ref{bootstrap CI}),  and show that for a properly chosen finite $B$, the nominal probability  coverage of size $1-\alpha$ will be achieved,  as $n$ approaches $+\infty$. We then use the result of the following
Theorem \ref{rate of bootstrap} to show that the rate at which the actual probability  coverage, of a classical bootstrap C.I. constructed using a finite number of  bootstrap sub-samples,  approaches its nominal coverage probability, $1-\alpha$, is no faster than  $n^{-1/2}$.

\begin{thm}\label{rate of bootstrap}
Assume that $E_{X}|X|^{3}<+\infty$. Consider a positive integer  $B\geq 2$ and let $\mathcal{H}$ be the class of all half space  subsets of $\mathbb{R}^{B}$. Define
\begin{equation*}
\mathbb{H}_n:=\Big( \frac{T_{n}(X-\mu)-T^{*}_{m_n}(1)}{\sqrt{2}}, \ldots, \frac{T_{n}(X-\mu)-T^{*}_{m_n}(B)}{\sqrt{2}} \Big)^{T},
 \end{equation*}
 and let $\mathbb{Y}:=(Z_1,\ldots,Z_B)^{T}$ be a $B$-dimensional Gaussian vector with mean $(0,\ldots,0)^{T}_{1\times B}$ and covariance matrix
\begin{equation}\label{covariance matrix}
\begin{pmatrix}
  1 & 1/2 & \cdots & 1/2 \\
  1/2 & 1 & \cdots & 1/2 \\
  \vdots  & \vdots  & \ddots & \vdots  \\
  1/2 & 1/2 & \cdots & 1
 \end{pmatrix}_{B\times B}.
\end{equation}
Then, as $n, \ m_n \rightarrow +\infty$ in such a way that    $m_n=o(n^2)$, for $\varepsilon>0$, the speed at which
$$
P_{w}\big\{\sup_{A\in \mathcal{H}} \Big| \bigotimes_{b=1}^{B}P_{X|w(b)}(\mathbb{H}_n \in A) -  P(\mathbb{Y}\in  A)  \Big|>\varepsilon   \big\}
$$
approaches zero is at best  $n^{-1/2}$, where $\bigotimes_{b=1}^{B}P_{X|w(b)}(.)$ is the conditional probability, given $\Big( \textbf{(}w_{1}^{(n)}(1),\ldots,w_{n}^{(n)}(1)\textbf{)},\ldots, \textbf{(}w_{1}^{(n)}(B),\ldots,w_{n}^{(n)}(B)\textbf{)}     \Big)$.
\end{thm}

\par
We note in passing that, on account of $m_n=o(n^2)$, Theorem \ref{rate of bootstrap} also holds true when $m_n=n$, i.e., when $B$ bootstrap sub-samples of  size $n$ are drawn from the original sample of size $n$.
Also, by virtue of Lemma 2.1 of S. Cs\"{o}rg\H{o} and Rosalsky \cite{Csorgo and Rosalsky}, the unconditional version of the conditional CLT of Theorem \ref{rate of bootstrap} continues to hold true under the same conditions as those of the  conditional one,  with the same rate of convergence that is    $n^{-1/2}$ at best.

\subsection{\normalsize{The classical bootstrap C.I.'s}}
The classical method of establishing  an asymptotic $1-\alpha$ size  bootstrap C.I. for $\mu$, as mentioned  before and formulated here,  is  based on the $B\geq 2$ ordered bootstrap readings $T^{*}_{n}[1]\leq \ldots \leq T^{*}_n [B]$, resulting from $B$ times independently re-sampling bootstrap sub-samples of size $m_n$, usually with $m_n =n$,    from the original sample with replacement by setting (cf., e.g., Efron and Tibshirani \cite{Efron and Tibshirani})
\begin{equation}\label{bootstrap CI}
T_{n}(X-\mu) \leq T^{*}_{ m_{n} }[\nu],
\end{equation}
where  $T^{*}_{m_{n}}[\nu]$ is the $\nu=(B+1).(1-\alpha)$th order statistic of  the $T^{*}_{m_{n}}(b)$, $1\leq b \leq B$, a bootstrap approximation to $\hat{\eta}_{n,\alpha}$ as in (\ref{(6.1)}). For simplicity,  we assume here that $\nu$ is an integer already.
\par
The so-called ideal bootstrap C.I. for $\mu$ is obtained when $B\rightarrow +\infty$. The validity of ideal bootstrap C.I.'s was established by Cs\"{o}rg\H{o} \emph{et al.} \cite{Scand} in terms of sub-samples of size $m_n$ when $E_{X}X^{2}<+\infty$ and it was also studied  previously by Hall \cite{Hall B} when $E_{X}|X|^{4+\delta}<+\infty$, with $\delta>0$. The common feature of  the results in the just mentioned two papers is that they both require that $n, m_n, B\rightarrow +\infty$, with $m_n=n$ in Hall \cite{Hall B}.
\par
In view of Theorem \ref{rate of bootstrap}, however, we establish another form of a bootstrap C.I. for $\mu$ of level $1-\alpha$ when $B$ is  fixed, as formulated  in our next section.
\par
\subsection{\normalsize{A  precise version of the classical  bootstrap C.I. for $\mu$ with fixed $B$ }}\label{Section 6.2}
We first consider  the counting  random variable $Y$  that counts the number of negative (or positive) components   $Z_t$, $1\leq t \leq B$,  in the Gaussian vector $\mathbb{Y}:=(Z_1,\ldots,Z_B)^T$. The distribution of $Y$ is
\begin{equation*}
P(Y\leq y)= \sum_{\ell=0}^{y} {B \choose \ell}  P(Z_{1}<0,\ldots,Z_{\ell}<0,Z_{\ell+1}>0,\ldots,Z_{B}>0 ),\
\end{equation*}
where, $y=0,\ldots,B$ and  $\mathbb{Y}:=(Z_1,\ldots,Z_B)^T$ has $B$-dimensional normal distribution with mean $(0,\ldots,0)^{T}_{1\times B} $ and covariance matrix as in (\ref{covariance matrix}).

\par
Let $y_{1-\alpha}$ be the $(1-\alpha)$th percentile of  $Y$, in other words,
\begin{equation*}
P(Y\leq y_{1-\alpha}-1)< 1-\alpha \ \textmd{and} \  \ P(Y\leq y_{1-\alpha})\geq 1-\alpha.
\end{equation*}
 We establish a more accurate  version of  $1-\alpha$  level classical bootstrap C.I. based on a finite number of,  $B$-times, re-sampling by setting
\begin{equation}\label{New bootstrap C.I.}
T_{n}(X-\mu)\leq T^{*}_{m_n}[y_{1-\alpha}],
\end{equation}
where,  $T^{*}_{m_n}[y_{1-\alpha}]$ is the $y_{1-\alpha}$th largest order statistic of the $B$ bootstrap versions  $T^{*}_{m_n}(1),\ldots,T^{*}_{m_n}(B)$ of $T^{*}_{m_n}$, constituting a new method for the bootstrap estimation of $\hat{\eta}_{n,\alpha}$ as in (\ref{(6.2)}), as compared to that of (\ref{bootstrap CI}).
\par
To show that as $n,m_n\rightarrow +\infty $, the probability coverage of the bootstrap C.I. (\ref{New bootstrap C.I.}) approaches its nominal probability of size $1-\alpha$, we first let
$$ (\bigotimes_{b=1}^{B}\Omega_{X,w(b)}, \bigotimes_{b=1}^{B}\mathfrak{F}_{X,w(b)},\bigotimes_{b=1}^{B}P_{X,w(b)}), \ 1\leq b \leq B,  $$
be the joint probability space of the $X$'s and the preceding array of the  weights $(w^{(n)}_{1}(b),\ldots,w^{(n)}_{n}(b))$, $1\leq b\leq B$, as in Theorem \ref{rate of bootstrap}, i.e.,   the  probability space generated by $B$ times, independently,  re-sampling from the original sample $X_1,\ldots,X_n$. Employing now the definition of order statistics, we  can compute  the  actual coverage probability of (\ref{New bootstrap C.I.}) as follows:
\\
\begin{eqnarray}
\alpha_{n}(B)&:=&\bigotimes_{b=1}^{B} P_{X,w(b)} \big\{  T_{n}\leq T^{*}_{m_n}[y_{1-\alpha}]    \big\}\nonumber\\
&=& \sum_{\ell=0}^{y_{1-\alpha}} {B\choose \ell} \bigotimes_{b=1}^{B} P_{X,w(b)} \Big\{ T_{n}(X-\mu) - T^{*}_{m_n}(1)\leq 0, \ldots, T_{n}(X-\mu) - T^{*}_{m_n}(\ell)\leq 0, \nonumber\\
&& ~~~~~~~~~~~~~~~~~~~~~~~~~~ T_{n}(X-\mu) - T^{*}_{m_n}(\ell+1)> 0,\ldots, T_{n}(X-\mu) - T^{*}_{m_n}(B)> 0 \Big\}\nonumber\\
&\longrightarrow& \sum_{\ell=0}^{y_{\alpha}} {B \choose \ell}  P(Z_{1}<0,\ldots,Z_{\ell}<0,Z_{\ell+1}>0,\ldots,Z_{B}>0 )=P(Y\leq y_{1-\alpha})  \nonumber
\end{eqnarray}
as\ $n,m_n\rightarrow +\infty$ such that $m_n=o(n^2)$.
\par
The preceding convergence  results from Theorem \ref{rate of bootstrap}  on assuming    the same conditions as those of  the latter theorem. It is noteworthy that the  conditional version of the preceding convergence,  in view of Theorem \ref{rate of bootstrap}, holds also true in probability-$P_w$ when one replaces the therein  joint probability  $\bigotimes_{b=1}^{B} P_{X,w(b)}$ by the conditional one $\bigotimes_{b=1}^{B}P_{X|w(b)}$.
\par
Recalling  that $y_{1-\alpha}$ is the $(1-\alpha)$th percentile of $Y$, we now conclude that \begin{equation}\label{useless benefit}
\alpha_{n}(B)\rightarrow \beta\geq 1-\alpha.
\end{equation}
This means that  the new bootstrap C.I. (\ref{New bootstrap C.I.}) is  valid.

\begin{remark}
The choice of the cutoff point $T^{*}_{m_n}[\nu]$    in the classical bootstrap C.I.
(\ref{bootstrap CI}) is done in a  blindfold way in comparison to the  more informed  choice of  $T^{*}_{m_n}[y_{1-\alpha}]$ in (\ref{New bootstrap C.I.}).
We also note that   increasing the value of $B$ to $B+1$ for the same confidence level $1-\alpha$,  only results in a different  cutoff point, which is a consequence of a change in the distribution of $Y$ which is now based on a $(B+1)$-dimensional Gaussian vector.
\end{remark}

\section{\normalsize{Numerical comparisons of the three C.I.'s for the population mean $\mu=E_{X} X$ }}\label{numerical comarison of the three}

The  numerical study below shows  that, in terms of the \emph{joint distribution} $P_{X,w}$, our confidence bound   of level 0.9000169 using the   bootstrapped pivot $G^{*}_{m_n}$, with $m_n=n$,  for $\mu$ outperforms the traditional confidence bound of level 0.9000169 with the pivot $T_{n}(X-\mu)$  and also  the classical bootstrapped confidence bound of the same level.
\par
We note that the classical bootstrap confidence bound,  as in  (\ref{bootstrap CI}),   coincides with our improved version of it (\ref{New bootstrap C.I.}) when $B=9$.

\par
In order to numerically compare  the performance of the three confidence bounds  for the population mean $\mu$ in terms of the joint distribution $P_{X,w}$, we let $m_n=n$ and  take a similar approach to the one used to illustrate Theorem \ref{Berry-Esseen} and its Corollary \ref{The Rate}. The only difference is that here we generate the weights $w_{i}^{(n)}$'s and the data  $X_1,\ldots,X_n$ \emph{simultaneously}.   More precisely, here we     generate 500 sets of the weights  $(w_{1}^{(n)}(b),\ldots,w_{n}^{(n)}(b))$ for  $b=1,\ldots,9$ and $X_1,\ldots,X_n$ at the same time.  We then compute the empirical distributions. The cycle is repeated 500 times. At the end,  we obtain the relative frequency of the empirical distributions that  did not deviate from the nominal limiting  probability $\Phi(1.281648)=0.9000169$ by more than 0.01. This procedure is  formulated as follows.

\begin{eqnarray*}\nonumber
{emp}_{X}T_{n}(X-\mu)&:=& \frac{\sum_{s=1}^{500} \mathds{1} \big\{ \big| 0.9000169-\frac{ \sum_{t=1}^{500} \mathds{1}(T^{(s,t)}_{n}(X-\mu)\leq 1.281648)}{500}\big|\leq 0.01 \big\}}{500}, \\
{emp}_{X,w}G^{*} &:=& \frac{\sum_{s=1}^{500} \mathds{1} \big\{\big| 0.9000169- \frac{\sum_{t=1}^{500} \mathds{1} (G^{*^{(s,t)}}_{m_n}\leq 1.281648)}{500}\big|\leq 0.01 \big\}}{500}, \\
{emp}_{X,w}Boot &:=& \frac{\sum_{s=1}^{500} \mathds{1}\big\{ \big| 0.9000169- \frac{\sum_{t=1}^{500} \mathds{1}
\big(T^{(s,t)}_{n}(X-\mu)\leq \max_{1\leq b \leq 9} T^{*^{(s,t)}}_{m_n}(b) \big)}{500}\big|\leq 0.01 \big\}}{500},
\end{eqnarray*}
where, for each $1\leq s\leq 500$, $T^{(s,t)}_{n}(X-\mu)$ and  $G^{*^{(s,t)}}_{m_n}$, with $m_n=n$,  stand for the respective values of $T_{n}(X-\mu)$ and  $G^{*}_{m_n}$, with $m_n=n$,  which are  computed, using 500 sets of $(X^{(s,t)}_1,\ldots,X^{(s,t)}_n)$ and $(w^{(n)^{(s,t)}}_1,\ldots,w^{(n)^{(s,t)}}_n)$, $1\leq t \leq 500$.

In a similar vein, for each $1\leq b \leq 9$ and  each $1\leq s \leq 500$, $T^{*^{(s,t)}}_{m_n}(b)$ represents the value of $T^{*}_{m_n}$ which are computed, using  the 500 simultaneously generated samples $(X^{(s,t)}_1,\ldots,X^{(s,t)}_n)$ and  $(w^{(n)^{(s,t)}}_{1}(b),\ldots,w^{(n)^{(s,t)}}_{n}(b))$ for  $ t =1, \ldots, 500$.

\par
The number 0.9000169 is the precise  nominal probability coverage of the  interval  $T_{n}(X-\mu)\leq \max_{1\leq b \leq 9} T^{*}_{m_n}(b)$ which, in view of our Theorem \ref{rate of bootstrap}, is $P(1\leq Y\leq 9)=1-P(Z_1>0,\ldots,Z_9>0)$, where  $Y$ is number of negative $Z_i$'s, $1\leq i \leq 9$,  in the $9$- dimensional Gaussian vector $(Z_1, \ldots,Z_9)^T$. To compute  the probability $P(Z_1>0,\ldots,Z_9>0)$,  we use  the    Genz algorithm  (cf. Genz \cite {Genz}), which is provided in the software R.

 \par
 It is noteworthy that the the examined classical bootstrap
  confidence bound for $\mu$, which is based on our modified version of the classical bootstrap confidence bound in (\ref{New bootstrap C.I.}) with $B=9$ bootstrap sub-samples, coincides with the  0.9 level classical  bootstrap confidence bound for $\mu$, of the form  (\ref{bootstrap CI}),  that was constructed by Hall \cite{Hall B}, based on the  same number of bootstrap sub-samples.

\par

We use the statistical software R to conduct the latter numerical study and present the results in the  following table.
\begin{table}[ht]\label{Table 2}
\caption{Comparing  the three confidence bounds  for $\mu$}
\begin{center}
\begin{tabular}{ c|c|c|c|c  }
\hline \hline
  Distribution &  $n$ & ${emp}_{X,w}G^{*}$ & ${emp}_{X}T_{n}(X-\mu)$ & ${emp}_{X,w}Boot$  \\
\hline
 \multirow{3}{*}{$Poisson(1)$}  &20& 0.48 & 0.302& 0.248 \\
                                &30& 0.494& 0.300& 0.33   \\
                                &40& 0.496& 0.350& 0.316  \\
                         \hline

 \multirow{3}{*}{$Lognormal(0,1)$}     &20& 0.028 & 0.000 & 0.000  \\
                                       &30& 0.048 & 0.000 & 0.004    \\
                                       &40& 0.058 & 0.000 & 0.002 \\
                         \hline

 \multirow{3}{*}{$Exponentia(1)$}   &20& 0.280 &0.026& 0.058   \\
                                    &30& 0.276 &0.026& 0.084  \\
                                    &40& 0.332 &0.048& 0.108   \\
\hline
\end{tabular}

\end{center}
\end{table}

\section{{Proofs}}\label{Proofs}

\subsection*{\normalsize{Proof of Theorem \ref{Berry-Esseen}}}
Due to similarity  of the two cases we only  give  the proof of part (A) of this theorem.
The proof  relies on the fact that via conditioning on the weights $w^{(n)}_i$'s, $\sum_{i=1}^{n} \big| \frac{w^{(n)}_i}{m_n}- \frac{1}{n}  \big| (X_i-\mu)$ as a sum of independent and non-identically distributed random variables. This in turn enables us to use a Berry-Esse\'{e}n type inequality for self-normalized sums of independents and non-identically distributed random variables. Also, some of the ideas in the proof are similar to those  of Slutsky's theorem.
\par
We now write
\begin{eqnarray}
G^{*}_{m_n}&=& \frac{\sum_{i=1}^{n}\big| \frac{w^{(n)}_i}{m_n}-\frac{1}{n}  \big| (X_i -\mu) }{\sigma \sqrt{\sum_{i=}^{n}} (\frac{w^{(n)}_i}{m_n}-\frac{1}{n})^{2}  } + \frac{\sum_{i=1}^{n}\big| \frac{w^{(n)}_i}{m_n}-\frac{1}{n}  \big| (X_i -\mu) }{\sigma \sqrt{\sum_{i=}^{n}} (\frac{w^{(n)}_i}{m_n}-\frac{1}{n})^{2}  } \big( \frac{\sigma}{S_n}-1 \big) \nonumber\\
&=:& Z_{m_n} + Y_{m_n}.
\end{eqnarray}
In view of the above setup, for $t\in \mathds{R}$ and $\varepsilon_1>0$, we have
\begin{eqnarray}
-P_{X|w} (|Y_{m_n}|>\varepsilon)&+&P_{X|w}(Z_{m_n} \leq t-\varepsilon)\nonumber\\
&\leq& P_{X|w}(G^{*}_{m_n} \leq t) \nonumber\\
&\leq & P_{X|w}(Z_{m_n} \leq t+\varepsilon)+ P_{X|w} (|Y_{m_n}|>\varepsilon).\label{eq 1 proofs}
\end{eqnarray}
Observe now that for $\varepsilon_2>0$ we have
\begin{equation}\label{eq 2 proofs}
P_{X|w} (|Y_{m_n}|>\varepsilon)
\leq P_{X|w} \big( |Z_{m_n}| > \frac{\varepsilon}{\varepsilon_1} \big)+ P_{X} \big(|S^{2}_n -\sigma^2|>\varepsilon_{1}^{2} \big).
\end{equation}
One can readily see that
\begin{eqnarray*}
P_{X|w} \big( |Z_{m_n}|> \frac{\varepsilon}{\varepsilon_1} \big) &\leq& (\frac{\varepsilon_2}{\varepsilon_1})^{2} \frac{\sum_{i=1}^{n} (\frac{w^{(n)}_i}{m_n}-\frac{1}{n})^{2} E_{X}(X_1 - \mu)^{2} }{\sigma^2 \sum_{i=1}^{n} (\frac{w^{(n)}_i}{m_n}-\frac{1}{n})^{2} }\\
&=& (\frac{\varepsilon_1}{\varepsilon})^{2}.
\end{eqnarray*}
Applying now the preceding in (\ref{eq 2 proofs}), in view of (\ref{eq 1 proofs}) can be replaced by
\begin{eqnarray}
&&- (\frac{\varepsilon_1}{\varepsilon})^{2} - P_{X} \big(|S^{2}_n -\sigma^2|>\varepsilon_{1}^{2} \big)+ P_{X|w} (Z_{m_n}\leq t-\varepsilon)\nonumber\\
&\leq& P_{X|w}(G^{*}_{m_n} \leq t) \nonumber\\
&\leq&(\frac{\varepsilon_1}{\varepsilon})^{2} + P_{X} \big(|S^{2}_n -\sigma^2|>\varepsilon_{1}^{2} \big)+P_{X|w} (Z_{m_n}\leq t+\varepsilon).\label{eq 3 proofs}
\end{eqnarray}
Now, the continuity of the normal distribution $\Phi$ allows us to choose  $\varepsilon_3>0$ so that  so that
$\Phi(t+\varepsilon)-\Phi(t)<\varepsilon_2$ and $\Phi(t-\varepsilon)-\Phi(t)>-\varepsilon_2$. This combined with (\ref{eq 3 proofs}) imply  that

\begin{eqnarray}
&&- (\frac{\varepsilon_1}{\varepsilon})^{2} - P_{X} \big(|S^{2}_n -\sigma^2|>\varepsilon_{1}^{2} \big)+ P_{X|w} (Z_{m_n}\leq t-\varepsilon)-\Phi(t-\varepsilon)-\varepsilon_2 \nonumber\\
&\leq& P_{X|w}(G^{*}_{m_n} \leq t)-\Phi(t)\nonumber \\
&\leq&(\frac{\varepsilon_1}{\varepsilon})^{2} + P_{X} \big(|S^{2}_n -\sigma^2|>\varepsilon_{1}^{2} \big)+P_{X|w} (Z_{m_n}\leq t+\varepsilon)-\Phi(t+\varepsilon)+\varepsilon_2.\nonumber\\
\label{eq 4 proofs}
\end{eqnarray}
We now use the Berry-Esse\'{e}en inequality for independent and not identically distributed random variables (cf., e.g., Serfling \cite{Berry-Esseen}) to write

\begin{equation*}
P_{X|w} (Z_{m_n}\leq t+\varepsilon_1)-\Phi(t+\varepsilon_1)\leq  \frac{  C E_{X}|X-\mu|^3}{\sigma^{3/2}}. \frac{\sum_{i=1}^n  | \frac{w^{(n)}_i}{m_n}-\frac{1}{n}  |^3 }{\big(  \sum_{i=1}^n  ( \frac{w^{(n)}_i}{m_n}-\frac{1}{n}  )^2   \big)^{3/2}}
\end{equation*}
and
\begin{equation*}
P_{X|w} (Z_{m_n}\leq t-\varepsilon_1)-\Phi(t-\varepsilon_1)\geq  \frac{ - C E_{X}|X-\mu|^3}{\sigma^{3/2}}. \frac{\sum_{i=1}^n  | \frac{w^{(n)}_i}{m_n}-\frac{1}{n}  |^3 }{\big(  \sum_{i=1}^n  ( \frac{w^{(n)}_i}{m_n}-\frac{1}{n}  )^2   \big)^{3/2}},
\end{equation*}
where $C$ is the universal constant of Berry-Esse\'{e}n inequality.
\par
Incorporating these approximations into  (\ref{eq 4 proofs}) we arrive at

\begin{eqnarray*}
&&\sup_{-\infty<t<+\infty} \big| P_{X|w}(G^{*}_{m_n} \leq t)-\Phi(t)  \big| \nonumber\\
&&\leq  (\frac{\varepsilon_1}{\varepsilon})^{2} + P_{X} \big(|S^{2}_n -\sigma^2|>\varepsilon_{1}^{2} \big)+ \frac{C E_{X}|X-\mu|^3}{\sigma^{3/2}}. \frac{\sum_{i=1}^n  | \frac{w^{(n)}_i}{m_n}-\frac{1}{n}  |^3 }{\big(  \sum_{i=1}^n  ( \frac{w^{(n)}_i}{m_n}-\frac{1}{n}  )^2   \big)^{3/2}}+\varepsilon_2.\nonumber
\end{eqnarray*}

From the preceding  relation  we conclude that

\begin{equation}\label{eq 5 proofs}
P_{w} \big( \sup_{-\infty<t<+\infty} \big| P_{X|w}(G^{*}_{m_n} \leq t)-\Phi(t)  \big|>\delta    \big) \leq P_{w} \big( \frac{\sum_{i=1}^n  | \frac{w^{(n)}_i}{m_n}-\frac{1}{n}  |^3 }{\big(  \sum_{i=1}^n  ( \frac{w^{(n)}_i}{m_n}-\frac{1}{n}  )^2   \big)^{3/2}} > \delta_n    \big)
\end{equation}
with $\delta_n$  as defined in the statement of  Theorem \ref{Berry-Esseen}.

\par
For $\varepsilon>0$, the right hand side of (\ref{eq 5 proofs}) is  bounded above by

\begin{eqnarray}
 && P_{w}\Big\{  \sum_{i=1}^{n} \big| \frac{w^{(n)}_i}{m_n}-\frac{1}{n}  \big|^3 >\frac{\delta_n(1-\varepsilon)^{\frac{3}{2}}(1-\frac{1}{n})^{\frac{3}{2} } }{m^{\frac{3}{2}}_n} \Big\}\nonumber\\
&&+ P_{w} \big(  \Big| \frac{m_{n}}{1-\frac{1}{n} } \sum_{i=1}^{n} \big( \frac{w^{(n)}_i}{m_n}-\frac{1}{n}  \big)^2 -1  \Big|>\varepsilon   \big)\nonumber\\
&&=: \Pi_1(n)+\Pi_2(n).\nonumber
\end{eqnarray}
We  bound $\Pi_{1}(n)$ above by
\begin{eqnarray}
&&\delta^{-2}_n(1-\varepsilon)^{-3}(1-\frac{1}{n})^{-3 } m^{-3 }_{n}  ( n  + n^2) E_{w} (w^{(n)}_1 -\frac{m_n}{n}    )^6 \nonumber\\
&&= \delta^{-2}_n(1-\varepsilon)^{-3}(1-\frac{1}{n})^{-3 } m^{-3 }_{n}  ( n  + n^2) \{ \frac{15 m^{3}_n}{n^3} + \frac{25 m^{2}_n}{n^2} + \frac{m_n}{n} \}. \label{eq 6 proofs}
\end{eqnarray}

As for $\Pi_{2}(n)$, recalling that $E_{w}\big(\sum_{i=1}^{n} (\frac{w^{(n)}_{i}}{m_n}-\frac{1}{n})^2 \big)=\frac{(1-\frac{1}{n})}{m_n}$,  an application of  Chebyshev's inequality yields
\begin{eqnarray}
\Pi_{2}(n)&\leq& \frac{m^{2}_n}{\varepsilon^{2} (1-\frac{1}{n})^2  } E_{w}\big( \sum_{i=1}^{n} (\frac{w^{(n)}_{i}}{m_n}-\frac{1}{n})^2 - \frac{(1-\frac{1}{n})}{m_n}   \big)^2 \nonumber\\
&=& \frac{m^{2}_n}{\varepsilon^{2} (1-\frac{1}{n})^2  }\ E_{w} \Big\{ \Big(\sum_{i=1}^n (\frac{w^{(n)}_i}{m_n}-\frac{1}{n})^2\Big)^2 - \frac{(1-\frac{1}{n})^2}{m^{2}_n} \Big\}^2  \nonumber\\
&=& \frac{m^{2}_n}{\varepsilon
^{2} (1-\frac{1}{n})^2  } \Big\{ n E_{w} \big( \frac{w^{(n)}_1}{m_n}-\frac{1}{n} \big)^4 + n(n-1) E_{w} \big[\big( \frac{w^{(n)}_1}{m_n}-\frac{1}{n} \big)^2 \big( \frac{w^{(n)}_2}{m_n}-\frac{1}{n} \big)^2\big] \nonumber\\
&-& \frac{(1-\frac{1}{n})^2}{m^{2}_n}   \Big\}.\label{for appendix 2}
\end{eqnarray}
We now use the fact  that $w^{(n)}$'s are multinomially distributed to compute  the preceding relation. After some algebra it turns out that it can be bounded above by
\begin{eqnarray}
&&\frac{m^{2}_n}{\varepsilon^{2} (1-\frac{1}{n})^2  } \Big\{ \frac{1-\frac{1}{n}}{n^3 m^{3}_n } + \frac{(1-\frac{1}{n})^4}{m^{3}_n}  + \frac{(m_n -1)(1-\frac{1}{n})^2}{n m^{3}_n} +  \frac{4(n-1)}{n^3 m_n}       +\frac{1}{m^{2}_n}\nonumber\\
&&~~~~~~~~~~~~~~~~~ - \frac{1}{n m^{2}_n} + \frac{n-1}{n^{3} m^{3}_{n}} + \frac{4(n-1)}{n^2 m^{3}_{n}} - \frac{(1-\frac{1}{n})^2}{m^{2}_n}        \Big\}.\label{eq 7 proofs}
\end{eqnarray}
Incorporating (\ref{eq 6 proofs}) and (\ref{eq 7 proofs}) into (\ref{eq 5 proofs}) completes the proof of
 part (A) of Theorem \ref{Berry-Esseen}. $\square$

\subsection*{\normalsize{Proof of Corollary \ref{The Rate}}}
The proof parts (A) and (B) of this corollary is an immediate consequence of Theorem \ref{Berry-Esseen}.
\par
To prove  parts (C) and (D) of this corollary, in view of Theorem \ref{Berry-Esseen}, it suffices to show that, for arbitrary  $\varepsilon_1, \varepsilon_2 >0$, as $n,m_n \to +\infty$,
\begin{equation}\label{Edmonton 1}
P_{w}\big( P_{X|w} ( |S^{*}_{m_n}-S^{2}_{n}|>\varepsilon_1)>\varepsilon_2  \big)=O( \frac{n}{m^{2}_n} ).
\end{equation}
To prove the preceding result we first note that
\begin{eqnarray*}
S^{*^2}_{m_n}-S^{2}_{n}&=&\sum_{1\leq i\neq j \leq n} \big(\frac{w^{(n)}_i w^{(n)}_j}{m_n (m_n -1)}-\frac{1}{n(n-1)} \big) \frac{(X_i-X_j )^2}{2}\\
&=& \sum_{1\leq i\neq j \leq n} \big(\frac{w^{(n)}_i w^{(n)}_j}{m_n (m_n -1)}-\frac{1}{n(n-1)} \big) \big(\frac{(X_i-X_j )^2}{2}-\sigma^2\big)
\end{eqnarray*}
By virtue of the preceding observation, we proceed with the proof of (\ref{Edmonton 1}) by first letting $d^{(n)}_{i,j}:=\frac{w^{(n)}_i w^{(n)}_j}{m_n (m_n -1)}-\frac{1}{n(n-1)}$ and writing

\begin{eqnarray}
&&P_{w}\big\{ P_{X|w} ( \big|  \sum_{1\leq i\neq j \leq n}  d^{(n)}_{i,j} \big(\frac{(X_i-X_j )^2}{2}-\sigma^2\big) \big| >\varepsilon_1)\varepsilon_2  \big\} \nonumber\\
&&\leq P_{w}\big\{ E_{X|w} \big(   \sum_{1\leq i\neq j \leq n}  d^{(n)}_{i,j} \big(\frac{(X_i-X_j )^2}{2}-\sigma^2\big)\big)^2 >\varepsilon^{2}_{1} \varepsilon_2 \big\}.\label{Edmonton 2}
\end{eqnarray}
Observe now that
\begin{eqnarray}
&&E_{X|w} \big(   \sum_{1\leq i\neq j \leq n}  d^{(n)}_{i,j} \big(\frac{(X_i-X_j )^2}{2}-\sigma^2\big)\big)^2 \nonumber \\
&=& E_{X}\big(\frac{(X_1-X_2 )^2}{2}-\sigma^2\big)^2  \sum_{1\leq i\neq j \leq n}  (d^{(n)}_{i,j})^2  \nonumber\\
&+& \sum_{\substack{1\leq i,j,k\leq n\\ i,j,k\ are \ distinct} } d^{(n)}_{i,j} d^{(n)}_{i,k}
 E_{X} \big( (\frac{(X_i-X_j)^2}{2}-\sigma^2)(\frac{(X_i-X_k)^2}{2}-\sigma^2)\big)\nonumber  \\
&+& \sum_{\substack{1\leq i,j,k,l\leq n\\ i,j,k,l\ are\ distinct} } d^{(n)}_{i,j} d^{(n)}_{k,l}
 E_{X} \big( (\frac{(X_i-X_j)^2}{2}-\sigma^2)(\frac{(X_k-X_l)^2}{2}-\sigma^2)\big).\nonumber\\
 \label{Edmonton 3}
\end{eqnarray}
We note that in the preceding relation, since $i,j,k$ are distinct, we have that
\begin{eqnarray*}
&&E_{X} \big( (\frac{(X_i-X_j)^2}{2}-\sigma^2)(\frac{(X_i-X_k)^2}{2}-\sigma^2)\big)\\
&=& E\big\{ E\big( \frac{(X_i-X_j)^2}{2}-\sigma^2|X_i\big) E\big( \frac{(X_i-X_k)^2}{2}-\sigma^2|X_i\big)  \big\}= \frac{E_{X} (X^{2}_{1}-\sigma^2)}{4}.
\end{eqnarray*}
Also, since $i,j,k,l$ are distinct, we have that

\begin{equation*}
E_{X} \big( (\frac{(X_i-X_j)^2}{2}-\sigma^2)(\frac{(X_k-X_l)^2}{2}-\sigma^2) \big)
=  E^{2}_{X} \big( \frac{(X_i-X_j)^2}{2}-\sigma^2)=0.
\end{equation*}
Therefore, in view of (\ref{Edmonton 3}) and (\ref{Edmonton 2}),  the proof of (\ref{Edmonton 1}) follows if we show that

\begin{equation}\label{Edmonton 4}
\sum_{1\leq i\neq j \leq n}  (d^{(n)}_{i,j})^2=O_{P_{w}}( \frac{1}{m^{2}_n})
\end{equation}
and
\begin{equation}\label{Edmonton 5}
\sum_{\substack{1\leq i,j,k\leq n\\ i,j,k\ are \ distinct} } d^{(n)}_{i,j} d^{(n)}_{i,k}=O_{P_{w}}(\frac{n}{m^{2}_n}).
\end{equation}
Noting that, as $n,m_n\to +\infty$,
\begin{equation*}
E_{w}\big\{ \sum_{1\leq i\neq j \leq n} (d^{(n)}_{i,j})^2\big\} \sim \frac{1}{m^{2}_n}
\end{equation*}
and

\begin{equation*}
E_{w} \big|  \sum_{\substack{1\leq i,j,k\leq n\\ i,j,k\ are \ distinct} } d^{(n)}_{i,j} d^{(n)}_{i,k} \big| \leq n^3 E_{w}(d^{(n)}_{1,2})^2 \sim \frac{n}{m^{2}_{n}}.
\end{equation*}
The  preceding two conclusions imply (\ref{Edmonton 4}) and (\ref{Edmonton 5}), respectively.
Now the proof of Corollary \ref{The Rate} is complete. $\square$

\subsection*{\normalsize{Proof of Theorem \ref{rate of bootstrap}}}
Once again  this proof  will be done by the means of conditioning on  the weights  $w^{(n)}$'s and this allows  us to think of them as constant coefficients for the $X$'s.
\par
When using the bootstrap to approximate the cutoff points, via repeated re-sampling, like in the  C.I. (\ref{New bootstrap C.I.}),
the   procedure can be described as  \emph{vectorizing} each centered observation. More precisely,  as a result  of drawing $B$ bootstrap sub-samples and each time  computing the value of $T^{*}_{m_n}$(b), $1\leq b \leq B$,  for each univariate random variable  $(X_i-\mu)$, $1\leq i \leq n$, we have  the following transformation.
\begin{equation*}
(X_i-\mu)\longmapsto  \left[
                    \begin{array}{c}
                      \frac{(w^{(n)}_i(1)/m_n-1/n)}{S_n\sqrt{\sum_{1\leq j\leq n} (w^{(n)}_{j}(1)/m_n-1/n)^2}}(X_i-\mu) \\
                      \vdots \\
                      \frac{(w^{(n)}_i(B)/m_n-1/n)}{S_n\sqrt{\sum_{1\leq j\leq n} (w^{(n)}_{j}(B)/m_n-1/n)^2}}(X_i-\mu) \\
                    \end{array}
                  \right]_{B\times 1}
\end{equation*}
Viewing the problem from this perspective and  replacing the sample variance $S_{n}^2$ by $\sigma^2$  result  in having
\begin{eqnarray*}
&&\left[
  \begin{array}{c}
   S_n  (T_{n}(X-\mu)-T^{*}_{m_n}(1))/( \sigma \sqrt{2}) \\
    \vdots \\
   S_n ( T_{n}(X-\mu)-T^{*}_{m_n}(B))/( \sigma \sqrt{2}) \\
  \end{array}
\right]_{B\times 1} \\
&=& \sum_{i=1}^n \left[
           \begin{array}{c}
             (  \frac{1}{\sigma \sqrt{n}}-  \frac{(w^{(n)}_i(1)/m_n-1/n)}{\sigma \sqrt{\sum_{1\leq j\leq n} (w^{(n)}_{j}(1)/m_n-1/n)^2}} ) \frac{(X_i-\mu)}{\sqrt{2}} \\
             \vdots \\
             (  \frac{1}{\sigma \sqrt{n}}-  \frac{(w^{(n)}_i(B)/m_n-1/n)}{\sigma \sqrt{\sum_{1\leq j\leq n} (w^{(n)}_{j}(B)/m_n-1/n)^2}} ) \frac{(X_i-\mu)}{\sqrt{2}} \\
           \end{array}
         \right]_{B\times 1} \\
&=: & \sum_{i=1}^n \left[
           \begin{array}{c}
             \big(  \frac{1}{ \sqrt{n}}-  a_{i,n}(1) \big) \frac{(X_i-\mu)}{\sigma \sqrt{2}} \\
             \vdots \\
             \big(  \frac{1}{ \sqrt{n}}-  a_{i,n} (B) \big) \frac{(X_i-\mu)}{\sigma  \sqrt{2}} \\
           \end{array}
         \right]_{B\times 1}
\end{eqnarray*}
Conditioning on $w^{(n)}_i$s, the preceding representation is viewed as  a  sum of $n$ independent but not identically distributed  $B$-dimensional random  vectors.  This, in turn, enables us to use Theorem  1.1 of Bentkus \cite{Bentkus} to derive the rate of the conditional  CLT  in Theorem \ref{rate of bootstrap}. For the sake of simplicity, we give the proof of this theorem only for $B=2$, as the  proof  essentially remains the  same for $B\geq 2$. Also, in the proof
we will consider  half spaces of the form $A=\{ (x,y)\in\mathbb{R}^2: x\leq a,\ y> b \}$, where $a,b \in \mathbb{R}$, as other forms of half spaces can be treated in the same vein. Moreover, in  what will follow we let $\| . \|$  stand for the Euclidean norm on $\mathbb{R}^2$.
\par
We now continue the proof by an application of Theorem 1.1 of Bentkus \cite{Bentkus} as follows.
\\
\begin{eqnarray}
&&\sup_{(x_1,x_2)\in \mathbb{R}^2} \Big| \bigotimes_{b=1}^{2}P_{X|w(b)} \left[
  \begin{array}{c}
   S_n  (T_{n}(X-\mu)-T^{*}_{m_n}(1))/( \sigma \sqrt{2}) \leq x_1 \\
   S_n ( T_{n}(X-\mu)-T^{*}_{m_n}(2))/( \sigma \sqrt{2}) > x_2 \\
  \end{array}
\right] \nonumber\\
&& \qquad \qquad \ \  - P\big( Z_1 \leq x_1, Z_2 >x_2 \big)  \Big| \nonumber\\
&&\leq c 2^{1/4} \sum_{i=1}^{n} E_{X|w} \big{\|} C_{n}^{-1/2}  \left[
           \begin{array}{c}
             \big(  \frac{1}{ \sqrt{n}}-  a_{i,n}(1) \big) \frac{(X_i-\mu)}{\sigma \sqrt{2}} \\
             \vdots \\
             \big(  \frac{1}{ \sqrt{n}}-  a_{i,n} (B) \big) \frac{(X_i-\mu)}{\sigma  \sqrt{2}} \\
           \end{array}
         \right]   \big{\|}^3,\label{new eq 1 proofs}
\end{eqnarray}
where $c$ is an absolute constant and
\begin{equation*}
C_{n}^{-1/2}= \left(
                \begin{array}{cc}
                  A_n & -B_n \\
                  -B_n & A_n \\
                \end{array}
              \right)
\end{equation*}
with
\begin{eqnarray*}
A_n&=& \frac{  1+ \sqrt{ 1-  A^{2}_{n}(1,2) }   }{ \sqrt{2  + 2 \sqrt{ 1-  A^{2}_{n}(1,2) }  }  \big(  2-A^{2}_{n}(1,2)+ \sqrt{ 1-  A_{n}^2(1,2) }\big)}
\\
B_n&=& \frac{ A^{2}_{n}(1,2)  \big(  1+ \sqrt{ 1-  A_{n}^{2}(1,2) } \  \big)}{ \sqrt{2  + 2 \sqrt{ 1-  A^{2}_{n}(1,2) }  }  \big(  2-A^{2}_{n}(1,2)+ \sqrt{ 1-  A_{n}^2(1,2) }\big)}.
\end{eqnarray*}
where $A_{n}(1,2):=  1/2 (\sum_{i=1}^n a_{i,n}(1) a_{i,n}(2)+1 )$.

We note that $C_{n}^{-1/2}$ is the inverse of a   positive square root of the covariance matrix of the vector \begin{equation*}
\Big(S_n  (T_{n}(X-\mu)-T^{*}_{m_n}(1))/( \sigma \sqrt{2}), S_n  (T_{n}(X-\mu)-T^{*}_{m_n}(2))/( \sigma \sqrt{2})  \Big),
\end{equation*}
which is
\begin{equation*}
\left(
  \begin{array}{cc}
    1 & A_{n}(1,2) \\
    A_{n}(1,2) & 1 \\
  \end{array}
\right).
\end{equation*}
Some algebra shows that the R.H.S. of (\ref{new eq 1 proofs}) is equal to
\begin{eqnarray*}
&&c 2^{-3/2} E_{X} |X_1-\mu|^3 \sum_{i=1}^n \Big\{  \frac{2(A_n -B_n)^2}{n}+ (A^{2}_n+ B^{2}_n) \big(a^{2}_{i,n}(1)+a^{2}_{i,n}(2) \big)\\
&& -2  \frac{(A_n -B_n)^2 \big(a_{i,n}(1)+a_{i,n}(2) \big)  }{\sqrt{n}} -4 A_n B_n    a_{i,n}(1) a_{i,n}(2)  \Big\}^{3/2}.
\end{eqnarray*}
In summary, so far,  we have shown that
\begin{eqnarray}
&& \sup_{(x_1,x_2)\in \mathbb{R}^2} \Big| \bigotimes_{b=1}^{2}P_{X|w(b)} \left[
  \begin{array}{c}
   S_n  (T_{n}(X-\mu)-T^{*}_{m_n}(1))/( \sigma \sqrt{2}) \leq x_1 \\
       S_n ( T_{n}(X-\mu)-T^{*}_{m_n}(2))/( \sigma \sqrt{2}) > x_2 \\
  \end{array}
\right]\nonumber\\
&& \qquad \qquad \ \ - P\big( Z_1 \leq x_1, Z_2 >x_2 \big)  \Big|\nonumber\\
&& \leq c 2^{-3/2} E_{X} |X_1-\mu|^3 \sum_{i=1}^n \Big\{  \frac{2(A_n -B_n)^2}{n}+ (A^{2}_n+ B^{2}_n) \big(a^{2}_{i,n}(1)+a^{2}_{i,n}(2) \big)\nonumber\\
&& -2  \frac{(A_n -B_n)^2 \big(a_{i,n}(1)+a_{i,n}(2) \big)  }{\sqrt{n}} -4 A_n B_n  \  a_{i,n}(1) a_{i,n}(2)  \Big\}^{3/2} \nonumber\\
&&=: R(n). \label{new eq 2 proofs} 
\end{eqnarray}
To investigate the speed at which $P_{w}(R(n)>\varepsilon)$ approaches zero, as $n,m_n \rightarrow +\infty$, we first note that there is no cancelation of terms in the general term of the sum in $R(n)$ for each fixed $n, m_n$. The other important observation concerns  $(A_n -B_n)^2$, which is the  coefficient of the term $1/n$ in $R(n)$. It can be shown that as $n,m_n \rightarrow +\infty$ in such a way that $m_n=o(n^2)$,
\begin{equation}\label{cf. Appenxix 2}
\sum_{i=1}^n a_{i,n}(1) a_{i,n} (2) \rightarrow 0 \ in \ probability- P_w
\end{equation}
(cf. Appendix 2 for details). The latter means    that, as $n,m_n \rightarrow +\infty$ in such a way that $m_n=o(n^2)$,  the following in probability-$P_w$ statement holds true.
\begin{eqnarray}
(A_n -B_n)^2 &\rightarrow& \Big( \big(  (1+\sqrt{3/4} )^2-1/2 (1+\sqrt{3/4})   \big)/ \big(  \sqrt{2+2\sqrt{3/4}}(7/4 + \sqrt{3}{4}) \big)    \Big)^2 \nonumber\\
&=:& D>0 \label{new eq 3 proofs}
\end{eqnarray}
The preceding shows that $A_n$ and $B_n$  do  not contribute to the speed at which $P_{w}(R(n)>\varepsilon)\rightarrow 0$. At this stage one can  see   that      $P_{w}(R(n)>\varepsilon)$ approaches zero at a rate no faster  than  $n^{-1/2}$. To further elaborate on the latter conclusion  we employ Markov's inequality followed by an application of Jensen's inequality  to write
\begin{eqnarray}
P_{w}(R(n)>\varepsilon) &\leq & \varepsilon^{-1} \ c\ 2^{-3/2} n \ E^{1/2}_{w} \Big| \frac{2(A_n -B_n)^2}{n}+ (A^{2}_n+ B^{2}_n) \big(a^{2}_{i,n}(1)+a^{2}_{i,n}(2) \big) \nonumber\\
&& -2  \frac{(A_n -B_n)^2 \big(a_{i,n}(1)+a_{i,n}(2) \big)  }{\sqrt{n}} -4 A_n B_n    a_{i,n}(1) a_{i,n}(2)  \Big|^{3} \nonumber\\
&\leq &  \varepsilon^{-1} \ c \ n^{-1/2} \ E^{1/2}_{w} \big| A_n - B_n  \big|^{3}+L(n,c,\varepsilon). \label{new eq 4 proofs}
\end{eqnarray}
It is easy to check that $A_n - B_n$ is uniformly bounded in $n$. By this, and in view of (\ref{new eq 3 proofs}),   the dominated convergence theorem implies that, as $n,m_n \rightarrow +\infty$,
$$E^{1/2}_{w} \big| A_n - B_n  \big|^{3} \rightarrow D^{3/2}.   $$
Now, it is clear that the R.H.S. of (\ref{new eq 4 proofs}) approaches zero
no faster than $n^{-1/2}$.
To complete the proof of this theorem, for $\varepsilon_1, \varepsilon_2>0$, we use a Slutsky type argument to arrive at the following  approximation.
\begin{eqnarray*}
&& -P_{X|w} \big( \big| \frac{T_{n}(X-\mu)-G_{m_n}^{*}(1)}{\sigma   \sqrt{2} } (\sigma/S_n -1)   \big|>\varepsilon_1     \big) \\
&&- P_{X|w} \big( \big| \frac{T_{n}(X-\mu)-G_{m_n}^{*}(2)}{\sigma   \sqrt{2} } (\sigma/S_n -1)   \big|>\varepsilon_2     \big)\nonumber\\
&&+ P_{X|w} \left(
              \begin{array}{c}
             \frac{T_{n}(X-\mu)-G_{m_n}^{*}(1)}{\sigma   \sqrt{2} }   \leq x_1 -\varepsilon_1 \\
               \frac{T_{n}(X-\mu)-G_{m_n}^{*}(2)}{\sigma   \sqrt{2} } > x_2 + \varepsilon_2 \\
              \end{array}
            \right) - \Phi(x_1 -\varepsilon_1, x_2 + \varepsilon_2 )\\
&& +     \Phi(x_1 -\varepsilon_1, x_2 + \varepsilon_2 ) - \Phi(x_1, x_2 )\\
&& \leq P_{X|w} \left(
              \begin{array}{c}
             \frac{T_{n}(X-\mu)-G_{m_n}^{*}(1)}{S_n   \sqrt{2} }   \leq x_1  \\
               \frac{T_{n}(X-\mu)-G_{m_n}^{*}(2)}{S_n   \sqrt{2} } > x_2 \\
              \end{array}
            \right)-\Phi(x_1,x_2)\\
&& \leq  P_{X|w} \big( \big| \frac{T_{n}(X-\mu)-G_{m_n}^{*}(1)}{\sigma   \sqrt{2} } (\sigma/S_n -1)   \big|>\varepsilon_1     \big) \\
&&+ P_{X|w} \big( \big| \frac{T_{n}(X-\mu)-G_{m_n}^{*}(2)}{\sigma   \sqrt{2} } (\sigma/S_n -1)   \big|>\varepsilon_2     \big)\nonumber\\
&&+ P_{X|w} \left(
              \begin{array}{c}
             \frac{T_{n}(X-\mu)-G_{m_n}^{*}(1)}{\sigma   \sqrt{2} }   \leq x_1 +\varepsilon_1 \\
               \frac{T_{n}(X-\mu)-G_{m_n}^{*}(2)}{\sigma   \sqrt{2} } > x_2 - \varepsilon_2 \\
              \end{array}
            \right) - \Phi(x_1 +\varepsilon_1, x_2 - \varepsilon_2 )\\
&&+ \Phi(x_1 +\varepsilon_1, x_2 - \varepsilon_2 ) - \Phi(x_1, x_2 ).
\end{eqnarray*}
By virtue of  (\ref{new eq 2 proofs}) and also by continuity of the  bivariate normal distribution function, for some $\varepsilon_3>0$,   we can replace the preceding approximations by

\begin{eqnarray}
&&\big| P_{X|w} \left(
              \begin{array}{c}
             \frac{T_{n}(X-\mu)-G_{m_n}^{*}(1)}{   \sqrt{2} }   \leq x_1  \\
               \frac{T_{n}(X-\mu)-G_{m_n}^{*}(2)}{   \sqrt{2} } > x_2 \\
              \end{array}
            \right)-\Phi(x_1,x_2) \big| \nonumber \\
&& \leq  P_{X|w} \big( \big| \frac{S_n\big(T_{n}(X-\mu)-G_{m_n}^{*}(1)\big)}{\sigma   \sqrt{2} } (\sigma/S_n -1)   \big|>\varepsilon_1     \big) \nonumber \\
 && + \quad P_{X|w} \big( \big| \frac{S_n\big(T_{n}(X-\mu)-G_{m_n}^{*}(2)\big)}{\sigma   \sqrt{2} } (\sigma/S_n -1)   \big|>\varepsilon_2     \big)\nonumber\\
&& + \quad R(n) +\varepsilon_3. \label{new eq 5 proofs}
\end{eqnarray}
One can show that, as $n,m_n \rightarrow +\infty$ so that $m_n=o(n^2)$,  $P_{X|w}(S_n(T_{n}(X-\mu)-G_{m_n}^{*})/\sigma   \sqrt{2} \leq t) \rightarrow \Phi(t)$ in probability-$P_{w}$ (cf. Appendix 3 for details). The latter implies that the first two terms in the R.H.S. of (\ref{new eq 5 proofs}) approach zero, as $n, m_n \to +\infty$. As we have already noted,   $R(n)$ goes to zero with a rate that as best is  $n^{-1/2}$ in probability-$P_w$. By sending $\varepsilon_1,\varepsilon_2 \rightarrow 0$,  $\varepsilon_3$ goes to zero too and this finishes  the proof of Theorem \ref{rate of bootstrap}. $\square$


\section{\normalsize{Appendix 1}}
Consider the original sample $(X_1,\ldots,X_n)$ and assume that the sample size $n\geq 1$ is fixed. It is known that the bootstrap estimator of the mean based on $B$ independent sub-samples $(X^{*}_{1}(b),\ldots,X^{*}_{n}(b))$, $1\leq b \leq B$ can be computed as
\begin{eqnarray}
\hat{X_{nB}^*}&=&\frac{\sum_{b=1}^{B} \bar{X^{*}}_n (b)}{B} \nonumber\\
&=& \frac{1}{n B} \sum_{i=1}^{n} \Big(\sum_{b=1}^{B} w^{(n)}_{i}(b) \Big) X_i, \label{Only one sub-sample 1}
\end{eqnarray}
where $w^{(n)}_{i}(b)$ is the $\#$ of times the index $i$, i.e., $X_i$ is chosen in the $b$th bootstrap sub-sample $(X^{*}_{1}(b),\ldots,X^{*}_{n}(b))$. Observe now that $\sum_{b=1}^{B} w^{(n)}_{i}(b)$  counts the total $\#$ of times $X_{i}$ has appeared in the $m:=n B$ bootstrap sub-samples. Also, observe that for fixed $n$,  $B\rightarrow +\infty$ is equivalent to $m\rightarrow +\infty$. Therefore, in view of  (\ref{Only one sub-sample 1}) we can write
\begin{equation}\label{Only one sub-sample 2}
\hat{X^{*}}_{nB}= \bar{X^{*}}_{m}.
\end{equation}
This means that taking a large number, $B$, of independent bootstrap sub-samples is equivalent to  taking only one large bootstrap sub-sample.
\par
We are now to show that when $n$ is fixed, as $m\rightarrow +\infty$, we have $\bar{X^{*}}_{m_n}\rightarrow \bar{X}_{n}$ in probability $P_{X,w}$. To do so, without loss of generality we assume that $\mu=0$, let $\varepsilon_1,\varepsilon_2>0$ and write
\begin{eqnarray}
P_{w} \big\{ P_{X|w} \big( \big|\bar{X^{*}}_{m}- \bar{X}_{n}\big|> \varepsilon_1 \big)>\varepsilon_2   \big\}
&\leq&    P_{w} \big\{ E_{X|w} \big( \sum_{i=1}^n (\frac{w^{(n)}_i}{m}-\frac{1}{n}) X_i \big)^2 >\varepsilon^{2}_1 \varepsilon_2     \big\}\nonumber\\
&=& P_{w} \big\{ \sum_{i=1}^n  (\frac{w^{(n)}_i}{m}-\frac{1}{n}) ^2 >  \sigma^{-2} \varepsilon^{2}_1 \varepsilon_2   \big\}\nonumber\\
&\leq& \sigma^{2} \varepsilon^{-2}_1 \varepsilon^{-1}_2 \ n E_{w} \big( \frac{w^{(n)}_1}{m}-\frac{1}{n} \big)^{2}\nonumber\\
&\leq& \sigma^{-2} \varepsilon^{-2}_1 \varepsilon^{-1}_2 \frac{(1-\frac{1}{n})}{m}\rightarrow 0, \ as\ m\rightarrow \infty. \nonumber\\
\label{(5.30)}
\end{eqnarray}
The preceding means that $P_{X|w} \big( \big|\bar{X^{*}}_{m_n}- \bar{X}_{n}\big|> \varepsilon_1 \big)\rightarrow 0$ in probability-$P_{w}$, hence from the dominated convergence theorem we conclude that $\bar{X^{*}}_{m}\rightarrow \bar{X}_{n}$ in probability $P_{X,w}$.  $\square$
\par
We are now to show that the bootstrap sample variance which we denote by $S^{*^2}_{m}$ is a in probability consistent estimator of the ordinary sample variance $S^{2}_{n}$ for each fixed $n$, when $m\rightarrow+\infty$. To do so, we denote the bootstrap sub-sample of size $m$ by $(X^{*}_1,\ldots,X^{*}_{m})$. By this we have that
$S^{*^2}_{m}=\sum_{k=1}^m \big( X^{*}_k-\bar{X}^{*}_{m} \big)^2\big/ (m-1)$. Employing now the $u$-statistic representation of the sample variance enables us to write
\begin{eqnarray*}
S^{*^2}_{m}&=& \frac{\sum_{1\leq k\leq l\leq m }  (X^{*}_{k}-X^{*}_{l})^{2} }{2 m(m-1)}\\
&=& \frac{\sum_{1\leq i\leq j\leq n } w^{(n)}_i w^{(n)}_j (X_{i}-X_{j})^{2} }{2 m(m-1)}.
\end{eqnarray*}
The preceding relation is the weighted form of $S^{*^2}_{m}$ and it is based on the fact that the terms $(X^{*}_{k}-X^{*}_{l})^2$, $1\leq k \neq l \leq m$, are replications  of $(X_i -X_j)^2$, $1\leq i \neq j \leq n$. Therefore, the deviation $S^{*^2}_{m}-S^{2}_n$ can be written as follows.
\begin{equation}\nonumber
 \sum_{1\leq i\neq j\leq n} \big( \frac{1}{2m(m-1)}-\frac{1}{2n(n-1)}\big) (X^{2}_i-X^{2}_j).
\end{equation}
Now for $\ep_1,\ep_2>0$ we write

\begin{eqnarray}
&&P_{w}\Big( P_{X|w}\big(\Big| S^{*^2}_{m}-S^{2}_n  \Big|>\ep_1\big)>\ep_2 \Big)  \nonumber\\
&&=P_w\big(  P_{X|w}\big( \big|\mathop{\sum\sum}_{1\leq i\neq j\leq n}
(\frac{w_i^{(n)}w_j^{(n)}}{m(m-1)}
- \frac{1}{ n(n-1)} ) (X_i-X_j)^2\big| > 2 \ep_1\big)>\ep_2  \big)\nonumber \\
&\leq&  P_w\big(   \mathop{\sum\sum}_{1\leq i\neq j\leq n}
\big| \frac{w_i^{(n)}w_j^{(n)}}{m_n(m_n-1)}
- \frac{1}{ n(n-1)} \big| E_{X}(X_i-X_j)^2 > 2 \ep_1 \ep_2  \big)\nonumber \\
&\leq&  P_w\big(   \mathop{\sum\sum}_{1\leq i\neq j\leq n}
\big| \frac{w_i^{(n)}w_j^{(n)}}{m_n(m_n-1)}
- \frac{1}{ n(n-1)} \big|  > \ep_1 \ep_2 \sigma^{-2}  \big).\label{(5.31)}
\end{eqnarray}

The preceding relation can be bounded above by:
\begin{eqnarray}
&& \ep^{-2}_1 \ep^{-2}_2 \sigma^{4}\Big\{ n(n-1) E_{w} \big(  \frac{w^{(n)}_1 w^{(n)}_2}{m (m-1)}-\frac{1}{n(n-1)} \big)^2 \nonumber \\
&+& n(n-1)(n-2) E_{w}\Big( \big|\frac{w^{(n)}_1 w^{(n)}_2}{m (m-1)}-\frac{1}{n(n-1)}\big|   \big|\frac{w^{(n)}_1 w^{(n)}_3}{m (m-1)}-\frac{1}{n(n-1)}\big|   \Big) \nonumber\\
&+& n(n-1)(n-2)(n-3) E_{w}\Big( \big|\frac{w^{(n)}_1 w^{(n)}_2}{m (m-1)}-\frac{1}{n(n-1)}\big|   \big|\frac{w^{(n)}_3 w^{(n)}_4}{m (m-1)}-\frac{1}{n(n-1)}\big|   \Big)  \Big\} \nonumber
\\
&\leq& \ep^{-2}_1 \ep^{-2}_2 \sigma^{4} \Big\{ n(n-1) E_{w} \big(  \frac{w^{(n)}_1 w^{(n)}_2}{m (m-1)}-\frac{1}{n(n-1)} \big)^2 \nonumber\\
&+& n(n-1)(n-2) E_{w} \big(  \frac{w^{(n)}_1 w^{(n)}_2}{m (m-1)}-\frac{1}{n(n-1)} \big)^2 \nonumber \\
&+& n(n-1)(n-2)(n-3)  E_{w} \big(  \frac{w^{(n)}_1 w^{(n)}_2}{m(m -1) } -\frac{1}{n(n-1)} \big)^2 \Big\} \nonumber \\
&=& \ep^{-2}_1 \ep^{-2}_2 \sigma^{4}\big\{ n(n-1)\nonumber\\
&+& n(n-1)(n-2)+ n(n-1)(n-2)(n-3) \big\}  \Big\{ \frac{1}{n^4 m^{2}  } + \frac{n}{n^4 m^{2} } + \frac{n^2}{ n^4 m^{2}} \Big\}.\nonumber
 \end{eqnarray}
Clearly, the preceding term approaches zero when $m\rightarrow +\infty$, for each fixed $n$. By this we have shown that $S^{*^2}_{m}\rightarrow S^{2}_{n}$ in probability-$P_{X,w}$, when $n$ is fixed and only $m\rightarrow +\infty$. $\square$
\subsection*{\normalsize{Consistency of $ \bar{X}_{n,m_n}^{*}$ in (\ref{(1.22)'})} }
We give the proof of  (\ref{(1.22)''}) for $m_n=n$, noting that the below proof remains the same for $m_n\leq n$ and it can  be  adjusted for the case $m_n= k n$, where $k$ is a positive integer.
In order establish (\ref{(1.22)''}) when $m_n=n$, we first note that
\begin{equation*}
E_{X|w}(\sum_{i=1}^n|\frac{w^{(n)}_i}{n}-\frac{1}{n}|)=2(1-\frac{1}{n})^{n}
\end{equation*}
and  for $\varepsilon_1, \varepsilon_2, \varepsilon_3>0$, we proceed  as follows.
\begin{eqnarray*}
&&P_{w} \big\{ P_{X|w} \big( \big| \bar{X}_{n,m_n}^{*}-\mu   \big|>\varepsilon_1        \big)>\varepsilon_2   \big\}\\
&\leq& P_{w}\big\{ P_{X|w} \Big( \big| \bar{X}_{n,m_n}^{*}-\mu   \big|>\varepsilon_1        \Big)>\varepsilon_2,  \big| \sum_{=1}^n|\frac{w^{(n)}_j}{n}-\frac{1}{n}|-2(1-\frac{1}{n})^{n}\big|\leq \varepsilon_3   \big\}\\
&+& P_{w} \big\{  \big| \sum_{j=1}^n|\frac{w^{(n)}_j}{n}-\frac{1}{n}|-2(1-\frac{1}{n})^{n}\big|>\varepsilon_3    \big\}\\
&\leq& P_{w} \big\{ P_{X|w} \big( \big|\sum_{i=1}^n |\frac{w^{(n)}_j}{n}-\frac{1}{n}|(X_i-\mu) \big|>  \varepsilon_1 \big(2(1-\frac{1}{n})^{n}-\varepsilon_3\big)          \big)>\varepsilon_2 \big\}\\
&+& \varepsilon_{3}^{-2}  E_{w} \big( \sum_{j=1}^n|\frac{w^{(n)}_j}{n}-\frac{1}{n}|-2(1-\frac{1}{n})^{n}   \big)^2\\
&\leq& P_{w} \big\{ \sum_{i=1}^n  (\frac{w^{(n)}_i}{n}-\frac{1}{n}) ^2 >  \sigma^{-2} \big(2(1-\frac{1}{n})^{n}-\varepsilon_3\big)^{2} \varepsilon_2   \big\}\\
&+& \varepsilon_{3}^{-2} \big\{  n E_{w} (\frac{w^{(n)}_1}{n}-\frac{1}{n}) ^2 +n(n-1) E_{w}\big( \big|  \frac{w^{(n)}_1}{n}-\frac{1}{n} \big|  \big| \frac{w^{(n)}_2}{n}-\frac{1}{n}  \big|  \big) - 4  (1-\frac{1}{n})^{2n}   \big\}\\
&=:& K_1(n)+K_2(n).
\end{eqnarray*}
A similar argument to that in (\ref{(5.30)}) implies that, as $n\to +\infty$, and then $\varepsilon_3\to 0$, $K_1(n) \to 0$. As for $K_2 (n)$, we note that
\begin{eqnarray*}
E_{w} (\frac{w^{(n)}_1}{n}-\frac{1}{n}) ^2&=&n^{-2} (1-\frac{1}{n})\\
E_{w}\big( \big|  \frac{w^{(n)}_1}{n}-\frac{1}{n} \big|  \big| \frac{w^{(n)}_2}{n}-\frac{1}{n}  \big|  \big)&=& -n^{-3}+4 n^{-2} (1-\frac{1}{n})^n (1-\frac{1}{n-1})^n.
\end{eqnarray*}
Observing now that, as $n\to +\infty$,
\begin{equation*}
n(n-1)E_{w}\big( \big|  \frac{w^{(n)}_1}{n}-\frac{1}{n} \big|  \big| \frac{w^{(n)}_2}{n}-\frac{1}{n}  \big|  \big)-4 (1-\frac{1}{n})^{2n}\to 0,
\end{equation*}
we imply that, as $n \to +\infty$, $K_2(n)\to 0$. By this we have concluded the consistency of  $\bar{X}_{n,m_n}^{*}$ for the population mean $\mu$, when $m_n=n$. $\square$

\section{\normalsize{Appendix 2}}

We are now to show that, as $n,m_n \rightarrow +\infty$ in such a way that $m_n/n^2 \rightarrow 0$, as  in Theorem \ref{rate of bootstrap}, (\ref{cf. Appenxix 2})  holds true. In order to do so, we let
  $\ep_1,\ep_2$ and $\ep_3$ be arbitrary positive numbers and  write
\begin{eqnarray*}
&&P\Big( \big| \sumn a_{i,n}(1) a_{i,n}(2) \big|>\ep_1    \Big)\\
&\leq& P \Big(  \frac{m_n}{(1-\frac{1}{n}) }\big| \sum_{i=1}^{n}  \big( \frac{w^{(n)}_{i}(1)}{m_n}-\frac{1}{n} \big) \big( \frac{w^{(n)}_{i}(2)}{m_n}-\frac{1}{n} \big)   \big|  >\ep_1(1-\ep_2)(1-\ep_3)    \Big)
\\
&+& P \Big(  \big| \frac{m_n}{(1-\frac{1}{n})} \sumn \big( \frac{w^{(n)}_{i}(1)}{m_n}-\frac{1}{n}  \big)^{2}-1  \big|>\ep_2  \Big) \qquad \qquad \\
&+& P \Big(  \big| \frac{m_n}{(1-\frac{1}{n})} \sumn \big( \frac{w^{(n)}_{i}(2)}{m_n}-\frac{1}{n}  \big)^{2}-1  \big|>\ep_3  \Big)\\
&=:& \pi_1 (n)+\pi_2 (n)+ \pi_3 (n).
\end{eqnarray*}
The last two terms in the preceding relation have already been shown to approach zero as $\frac{m_n}{n^{2}}\to 0$ (cf. (\ref{for appendix 2})). We now show that the first term approaches zero as well  in view of the following argument which relies on the facts that  $w^{(n)}_{i}$'s, $1\leq i \leq n$ are multinoialy distributed and that for each $1\leq i,j \leq n$, $w^{(n)}_{i}(1)$  and $w^{(n)}_{j}(2)$ are i.i.d.  (in terms of $P_w$).
\par
To show that  $\pi_1 (n)=o(1)$, as $n,m_n \rightarrow +\infty$ such that $m_n/n^2 \rightarrow 0$,  in what will follow,  we put  $\varepsilon_4:= \varepsilon_1(1-\varepsilon_2)(1-\varepsilon_3)$ to write
\begin{eqnarray*}
\pi_1 (n)&\leq& \varepsilon^{-2}_4 \frac{m^{2}_{n}}{(1-\frac{1}{n})^2} \Big\{   n E^{2}_{w} \big( \frac{w^{(n)}_{1}(1)}{m_n}-\frac{1}{n} \big)^{2} + \\
&&  n(n-1) E^{2}_{w} \Big[ \big(\frac{w^{(n)}_{1}(1)}{m_n}-\frac{1}{n}\big) \big(\frac{w^{(n)}_{2}(1)}{m_n}-\frac{1}{n}\big)  \Big]      \Big\}\\
&=& \varepsilon^{-2}_4 \frac{m^{2}_{n}}{(1-\frac{1}{n})^2} \Big\{   n\big( \frac{(1-\frac{1}{n})}{n m_n} \big)^{2}+
n(n-1) \big(  \frac{-1}{m_n n^{2}} \big)^{2}    \Big\}\\
&\leq& \varepsilon^{-2}_4 \Big( \frac{1}{n}+\frac{1}{n^{2}(1-\frac{1}{n})^2  } \Big)\to 0.
\end{eqnarray*}

\par
The preceding result completes the proof of (\ref{cf. Appenxix 2}). $\square$

\section{\normalsize{Appendix 3}}
Noting that the expression $S_n(T_{n}(X-\mu)-G_{m_n}^{*})/\sigma   \sqrt{2} $ can be written as
$$ S_n(T_{n}(X-\mu)-G_{m_n}^{*})/\sigma   \sqrt{2} = \sum_{i=1}^n  ( \frac{1}{\sqrt{n}} - \frac{(  \frac{w^{(n)}_i}{m_n} -\frac{1}{n})}{\sqrt{\sum_{j=1}^n   (\frac{w^{(n)}_i}{m_n} -\frac{1}{n})^2   } } )  (\frac{X_i - \mu}{\sigma \sqrt{2}} ) $$
makes it clear that, in view of Lindeberge-Feller CLT
, in order to have $P_{X|w}(S_n(T_{n}(X-\mu)-G_{m_n}^{*})/\sigma   \sqrt{2}\leq t) \to \Phi(t)$ in probability-$P_w$, it suffices to show that, as  $n,m_n \to +\infty$ so that $m_n=o(n^2)$,
\begin{equation}\label{new eq 6 proofs}
\max_{1\leq i \leq n}  \Big| \frac{(  \frac{w^{(n)}_i}{m_n} -\frac{1}{n})}{\sqrt{\sum_{j=1}^n   (\frac{w^{(n)}_i}{m_n} -\frac{1}{n})^2   } } - \frac{1}{\sqrt{n}}  \Big|=o_{P_{w}}(1).
\end{equation}
Considering that the L.H.S. of (\ref{new eq 6 proofs}) is bounded above by
\begin{equation*}
\max_{1\leq i \leq n}   \frac{  |  \frac{w^{(n)}_i}{m_n} -\frac{1}{n} |}{\sqrt{\sum_{j=1}^n   (\frac{w^{(n)}_i}{m_n} -\frac{1}{n})^2   } } + \frac{1}{\sqrt{n}} ,
\end{equation*}
(\ref{new eq 6 proofs}) follows if one shows that, as $n,m_n \rightarrow +\infty$ so that $m_n=o(n^2)$,

\begin{equation}
\max_{1\leq i \leq n}   \frac{  |  \frac{w^{(n)}_i}{m_n} -\frac{1}{n} |}{\sqrt{\sum_{j=1}^n   (\frac{w^{(n)}_i}{m_n} -\frac{1}{n})^2   } }=o_{P_{w}}(1).
\end{equation}
In order to establish the latter,  for $\varepsilon, \varepsilon^{\prime} > 0$, we write:
\begin{eqnarray*}
&&P_{w}\big(    \frac{\max_{1\leq i\leq n} \big( \frac{w_i^{(n)}}{m_n} - \frac{1}{n} \big)^2}
       {\sum_{i=1}^{n} \big( \frac{w_i^{(n)}}{m_n} - \frac{1}{n} \big)^2}>\ep  \big)\\
&\leq&    P_{w}\big(    \frac{\max_{1\leq i\leq n} \big( \frac{w_i^{(n)}}{m_n} - \frac{1}{n} \big)^2}
       {\sum_{i=1}^{n} \big( \frac{w_i^{(n)}}{m_n} - \frac{1}{n} \big)^2}>\ep, \big| \frac{m_n}{(1-\frac{1}{n})} \sum_{i=1}^{n} \big( \frac{w_i^{(n)}}{m_n} - \frac{1}{n} \big)^2 -1 \big|\leq {\ep}^{\prime}  \big)\\
&+& P_{w}\big(    \big| \frac{m_n}{(1-\frac{1}{n})} \sum_{i=1}^{n} \big( \frac{w_i^{(n)}}{m_n} - \frac{1}{n} \big)^2 -1 \big|> {\ep}^{\prime}  \big)\\
\end{eqnarray*}

\begin{eqnarray*}
&=& P_w \big( \max_{1\leq i\leq n} \big( \frac{w_i^{(n)}}{m_n} - \frac{1}{n} \big)^2>\frac{\ep(1-\ep^{\prime})(1-\frac{1}{n})}{m_{n}}   \big)
\\
&+& P_w \big(    \big|  \sum_{i=1}^{n} \big( \frac{w_i^{(n)}}{m_n} - \frac{1}{n} \big)^2 -\frac{(1-\frac{1}{n})}{m_n}\big|> \frac{{\ep}^{\prime} (1-\frac{1}{n})   }{m_n}    \big).
\end{eqnarray*}
An upper bound for $P_w \big( \max_{1\leq i\leq n} \big( \frac{w_i^{(n)}}{m_n} - \frac{1}{n} \big)^2>\frac{\ep(1-\ep^{\prime})(1-\frac{1}{n})}{m_{n}}   \big)$ is
\begin{eqnarray*}
&& n   P_w  \big( \big| \frac{w_1^{(n)}}{m_n} - \frac{1}{n}  \big|>\sqrt{\frac{\ep(1-\ep^{\prime})(1-\frac{1}{n})}{m_{n}} }   \big)
\\
&\leq& n \exp\{- \sqrt{m_n}\ . \ \frac{\ep(1-\ep^{\prime})  (1-\frac{1}{n})  }{2\big( \frac{\sqrt{m_n}}{n}+ \sqrt{\ep(1-\ep^{\prime})  (1-\frac{1}{n})  } \big)}  \}.
\end{eqnarray*}
The preceding relation, which  is due to Bernstien's inequality, is a general term of a finite series when $m_n=o(n^2)$.
\par
$P_w \big(    \big|  \sum_{i=1}^{n} \big( \frac{w_i^{(n)}}{m_n} - \frac{1}{n} \big)^2 -\frac{(1-\frac{1}{n})}{m_n}\big|> \frac{{\ep}^{\prime} (1-\frac{1}{n})   }{m_n}    \big)$ was already shown, (cf. (\ref{for appendix 2})),   to approach zero as $n,m_n \to +\infty$ such that $m_n=o(n^2)$.

\section{\normalsize{Appendix 4}}
Viewing the bootstrapped mean  as a randomly weighted partial sum, allows one to think  about randomly weighted partial sums of the general form
\begin{equation}\label{Baysian bootstrap}
\sum_{i=1}^n v_{i}^{(n)} X_{i},
\end{equation}
where, $X_i$'s are i.i.d. random variables and $v_{i}^{(n)}$'s  are random weights which are independent from $X_i$'s and posses the properties that  $E_v (v_{i}^{(n)}/m_{n})=1/n$ and $\sum_{i=1}^n v_{i}^{(n)}=m_n$.  The motivation behind studying the latter sums is that, in addition to bootstrapping,  they allows considering the problem of stochastically re-weighing (designing) the observations using a  random weights $v_{i}^{(n)}$ whose distribution is usually known. Naturally, in this setup, on taking  $v_{i}^{(n)}:=w_{i}^{(n)}$ for each $1\leq i \leq n$, $n\geq 1$, i.e., on assuming the random weights to be multinomially distributed, the randomly weighted partial sum in (\ref{Baysian bootstrap}) coincides with the  bootstrapped mean as in $(\ref{bootstrap mean})$. This idea  first appeared in the area of  Baysian bootstrap (cf. for example Rubin  \cite{Rubin}).  The convergence in distribution of the partial sums of  (\ref{Baysian bootstrap}) were also studied by Cs\"{o}rg\H{o} \emph{et al.} \cite{Scand} via conditioning on the weights  (cf. Theorem 2.1 and Corollary 2.2 therein).
Noting that in the latter results only randomly weighted statistics similar to $T_{m_n}^{*}$ which are natural pivots for the sample variance $\bar{X}_n$  were studied. In view of the fact that $G_{m_n}^{*}$, as defined by (\ref{G^{*}}), was seen to be a natural pivot for the population mean $\mu:=E_X X$, in a similar fashion and to Theorem 2.1 and its Corollary 2.2 of  Cs\"{o}rg\H{o} \emph{et al.} \cite{Scand}, here we state a more general  conditional CLT, given the weights $v_{i}^{(n)}$'s, for the partial sums $\sum_{i=1}^n |\frac{v^{(n)}_{i}}{m_n} - \frac{1}{n}|(X_i-\mu)$, noting that the proofs of these results almost identical to those of Theorem 2.1 and its Corollary 2.2 of  Cs\"{o}rg\H{o} \emph{et al.} \cite{Scand} in view of the more general setup on notations in the latter paper. In order to estate the results we first generalize the definition of $G_{m_n}^{*}$ and $G_{m_n}^{**}$, as defined in (\ref{G^{*}}) and (\ref{G^{**}}), respectively,  to become
\begin{equation}\label{generalizing G^*}
\mathcal{G}^{*}_{m_n}:= \frac{\sum_{i=1}^{n} \big| \frac{v^{(n)}_{i}}{m_{n}} -\frac{1}{n}   \big| \big( X_{i}-\mu \big)   }{S_{n}\sqrt{\sum_{i=1}^{n} (\frac{v^{(n)}_i}{m_n}-\frac{1}{n})^{2}}  },
\end{equation}

\begin{equation}\label{generalizing G^**}
\mathcal{G}_{m_n}^{**}:=\frac{\sum_{i=1}^n |\frac{v_{i}^{(n)}}{m_n}-\frac{1}{n}|(X_i-\mu)}{S_{m_n}^{*} \sqrt{\sum_{i=1}^{n} (\frac{v^{(n)}_i}{m_n}-\frac{1}{n})^{2}}  },
\end{equation}
where, $m_n=\sum_{i=1}^n v_{i}^{(n)}$.

\begin{thm}\label{thm1}
Let $X, X_{1},X_{2},\ldots$ be real valued i.i.d. random variables with mean $0$ and variance   $\sigma^2,$  and assume that  $0<\sigma^2 < \infty$. Put $V_{i,n}:=  \big| \big( \frac{v^{(n)}_{i}}{m_n}-\frac{1}{n} \big) X_i \big|, \ 1\leq i \leq n$,  $V^{2}_{n}:= \sum_{i=1}^{n}\big( \frac{v^{(n)}_{i}}{m_n}-\frac{1}{n} \big)^{2}$, $M_n:=\frac{ \max_{1\leq i\leq n} \big( \frac{v_i^{(n)}}{m_{n}} - \frac{1}{n}\big)^2 }{\sumn \big( \frac{v_i^{(n)}}{m_{n}} - \frac{1}{n} \big)^2}$,  and let $Z$ be a standard normal random variable throughout. Then as, $n,m_n \to \infty$, having
\\
\begin{eqnarray}
&& M_n=o(1)\ a.s.-P_v   \label{A1}   \\
&&is\ equivalent\ to\ concluding\ the\ respective\ statements\ of\nonumber\\
 &&(\ref{A2})\ and \ (\ref{A3})\nonumber \ simultaneously\ as \ follows \nonumber\\
&& P_{X|v} \left(\mathcal{G}^*_{m_n} \leq t\right) \longrightarrow
P(Z\leq t)\ a.s.- P_{v}\ for \ all\  t \in \mathds{R} \label{A2}\\
&&and \nonumber \\
&&\max_{1\leq i \leq n} P_{X|v}(V_{i,n}\big/ (S_{n} V_{n}) >\ep)=o(1)\ a.s.-P_{v}, for \ all\  \ep>0, \label{A3}
\end{eqnarray}
and, in a similar vein, having
\begin{eqnarray}
&& M_n=o_{P_v}(1)  \label{B1}  \\
&& is\ equivalent\ to\ concluding\ the\ respective\ statements\ of\nonumber\\
&& (\ref{B2})\ and \ (\ref{B3})\  as\ below\ simultaneously \nonumber\\
&& P_{X|v} \left(\mathcal{G}^*_{m_n} \leq t\right) \longrightarrow
P(Z\leq t)\ in\ probability - P_{v}\ for \ all\  t \in \mathds{R} \nonumber\\
\label{B2}\\
&&and \nonumber\\
&&\max_{1\leq i \leq n} P_{X|v}(V_{i,n}\big/ (S_{n} V_{n}) >\ep)=o_{P_v}(1),for \ all\  \ep>0.  \label{B3}
\end{eqnarray}
Moreover,    assume that, as $n, m_n \to \infty$, we have for\ any\ $\ep>0$,
\begin{numcases}{
\ P_{X|v}\Big(\big|\frac{S_{m_n}^{*2}\Big/ m_n}{ \sigma^2\ \sumn\big( \frac{v_i^{(n)}}{m_n}-\frac{1}{n} \big)^2} -1 \big|>\varepsilon \Big)=}
o(1)\ a.s.-P_{v} \label{E}& \\
o_{P_{v}}(1). \label{F}&
\end{numcases}
Then, as  $n, m_n \to \infty$,  via (\ref{E}),   the  statement of (\ref{A1}) is also equivalent to having (\ref{c2}) and (\ref{c3}) simultaneously as below

\begin{eqnarray}
&& P_{X|v} \left(\mathcal{G}^{**}_{m_n} \leq t\right) \longrightarrow
P(Z\leq t)\ a.s.- P_{v}\ for \ all\  t \in \mathds{R} \label{c2}    \\
&& and \nonumber\\
&&\max_{1\leq i \leq n} P_{X|v}(V_{i,n}\big/ (S^{*}_{m_n}/\sqrt{m_n} )>\ep)=o(1)\ a.s.-P_{v}, \ for\ all\ \ep>0, \nonumber\\
 \label{c3}
\end{eqnarray}
and, in a similar  vein, via (\ref{F}), the statement (\ref{B1}) is also equivalent to having  (\ref{d2}) and (\ref{d3}) simultaneously as below

\begin{eqnarray}
&& P_{X|v} \left(T^{**}_{m_n} \leq t\right) \longrightarrow
P(Z\leq t)\ in\ probability - P_{v}\ for \ all\  t \in \mathds{R} \nonumber\\
\label{d2} \\
&&and   \nonumber \\
&&\max_{1\leq i \leq n} P_{X|v}(V_{i,n}\big/ (S^{*}_{m_n}/\sqrt{m_n} )>\ep)=o_{P_v}(1),\ for \ all \ \ep>0.  \label{d3}
\end{eqnarray}
\end{thm}
For verifying the technical conditions (\ref{E}) and (\ref{F})  as above, one does not need to know the actual finite   value  of $\sigma^2$.

\par
Now suppose that  $\displaystyle{v^{(n)}_{i}=\zeta_{i} }$, $1\leq i \leq n$, where $\zeta_{i}$ are positive i.i.d. random variables. In this case, noting that  the bootstrapped $t$-statistic $\mathcal{G}^{*}_{m_n}$ defined by (\ref{generalizing G^*}) is of the form:
\begin{equation}
\mathcal{G}^*_{m_n} =
\frac{ \displaystyle \sum^n_{i=1} | \frac{\zeta_{i}}{m_n} - \frac{1}{n} \big| (X_i-\mu)}
{S_n \sqrt{   \displaystyle\mathop{\sum}_{i=1}^{n} (\frac{\zeta_{i}}{m_n}-\frac{1}{n})^2 } },\label{T*zeta}
\end{equation}
where $\displaystyle{m_n=\sum_{i=1}^{n} \zeta_i}$.
\par
The following Corollary \ref{cor2} to  Theorem \ref{thm1}  establishes the validity of this scheme of bootstrap for $\mathcal{G}^{*}_{m_n}$, as defined by (\ref{T*zeta}), via conditioning on the   weights $\zeta_i$'s.

\begin{corollary}\label{cor2}
Assume that $0 < \sigma^{2}= var(X) <\infty$, and let $\zeta_1,\zeta_2,\ldots$ be a sequence of positive i.i.d. random variables  which are independent of $X_{1},X_{2}, \dots$ .  Then,    as $ n\to\infty$,
\\
(a) if   $E_{\zeta}(\zeta^{4}_{1})< \infty$,  then, mutatis mutandis, (\ref{A1}) is equivalent to having (\ref{A2}) and (\ref{A3}) simultaneously and, spelling  out only  (\ref{A3}),   in this context it reads
\begin{equation}
P_{X|\zeta} (\mathcal{G}_{m_n}^* \leq t) \l
ongrightarrow P(Z\leq t) \ \hbox {a.s.}-P_\zeta, \ for\ all\ t \in \mathds{R},
\end{equation}
(b) if   $E_{\zeta}(\zeta^{2}_{1})< \infty$,  then, mutatis mutandis, (\ref{B1}) is equivalent (\ref{B2}) and (\ref{B3})  simultaneously, and spelling  out only  (\ref{B2}),   in this context it reads
\begin{equation}
P_{X|\zeta} (\mathcal{G}_{m_n}^* \leq t) \longrightarrow P(Z\leq t) \ in \ probability-P_\zeta, \ for\ all\ t \in \mathds{R},
\end{equation}
where $Z$ is a standard normal random variable.
\end{corollary}


\end{document}